\documentclass[twoside,11pt]{amsart}

\usepackage{amsmath,latexsym,amssymb, times, enumerate}
\input amssym.def
\input amssym
\input xypic
\input xy
\xyoption{all}
\def\demo{\noindent{\bf Proof. }}
\def\sqr#1#2{{\vcenter{\hrule height.#2pt
        \hbox{\vrule width.#2pt height#1pt \kern#1pt
                \vrule width.#2pt}
        \hrule height.#2pt}}}
\def\square{\mathchoice\sqr64\sqr64\sqr{4}3\sqr{3}3}
\def\QED{\hfill$\square$}

\def\tratto{\mbox{\rule{2mm}{.2mm}$\;\!$}}

\def\m{{\mathfrak m}}
\def\n{{\mathfrak n}}

\def\p{{\mathfrak p}}

\newtheorem{Theorem}{Theorem}[section]
\newtheorem{Claim}{Claim}

\newtheorem{Lemma}[Theorem]{Lemma}
\newtheorem{Corollary}[Theorem]{Corollary}
\newtheorem{Proposition}[Theorem]{Proposition}

\newtheorem{Discussion}[Theorem]{Discussion}
\newtheorem{Notation and Discussion}[Theorem]{Notation and Discussion}
\newtheorem{Assumptions and Discussion}[Theorem]{Assumptions and Discussion}

\newtheorem{Example}[Theorem]{Example}
\newtheorem{Definition}[Theorem]{Definition}

\begin{document}

\baselineskip=16pt

\title[Generalized stretched ideals and  Sally's Conjecture]
{\Large\bf Generalized stretched ideals and Sally's Conjecture}

\author[P. Mantero  and Y. Xie]
{Paolo Mantero   \and Yu Xie}

\thanks{AMS 2010 {\em Mathematics Subject Classification}.
Primary 13A30; Secondary 13H15, 13B22, 13C14, 13C15, 13C40.}


\address{Department of Mathematics, Purdue University,
West Lafayette, Indiana 47907} \email{pmantero@math.purdue.edu}

\address{Department of Mathematics, University of Notre Dame,
Notre Dame, Indiana 46556} \email{yxie@nd.edu}

\vspace{-0.1in}

\begin{abstract}
Given a finite module $M$ over a Noetherian local ring $(R, \m)$,  we introduce the concept of $j$-stretched ideals on $M$. Thanks to a crucial specialization lemma, we show that this notion greatly generalizes (to arbitrary ideals, and with respect to modules) the classical definition of stretched $\m$-primary ideals of Sally and Rossi-Valla, as well as the notion of minimal and almost minimal $j$-multiplicity given  recently by Polini-Xie.
For $j$-stretched ideals $I$ on a Cohen-Macaulay module $M$, we show that ${\rm gr}_I(M)$ is Cohen-Macaulay if and only if two classical invariants of $I$, the reduction number and the index of nilpotency, are equal.
Moreover, for the same class of ideals, we provide a generalized version of Sally's conjecture (proving the almost Cohen-Macaulayness of  associated graded rings).
Our work unifies the approaches of Rossi-Valla and  Polini-Xie   and
generalizes  simultaneously results on the (almost) Cohen-Macaulayness
of associated graded  modules by several authors, including Sally, Rossi-Valla, Wang, Elias, Rossi, Corso-Polini-Vaz Pinto, Huckaba and Polini-Xie.
\end{abstract}

\maketitle

\vspace{-0.2in}

\section{Introduction}

The maximal ideal $\n$ of an Artinian local ring $A$ is called {\em stretched} if $\n^2$ can be generated by one element, i.e., if $\n^2$ is a principal ideal. The maximal ideal $\m$ of a Cohen-Macaulay local ring $R$ is dubbed {\em stretched} if the image of $\m$ in an Artinian reduction of $R$ is stretched.
The notion of stretchedness was introduced by Sally in \cite{S3} to generalize Cohen-Macaulay local rings of minimal or almost minimal multiplicity
, two classes of rings for which she proved results on the  Cohen-Macaulayness of the associated graded rings of the maximal ideals.
In fact, in the case of minimal multiplicity, she proved that  the associated graded ring of the maximal ideal is always Cohen-Macaulay \cite{S1}. She then extended this result by proving a clean criterion for the Cohen-Macaulayness of the associated graded rings of  stretched maximal ideals \cite{S3}. Furthermore, she  proved that in the case of almost minimal multiplicity, the associated graded ring is Cohen-Macaulay, provided that the type of the ring is not maximal \cite{S4}. However, she also found examples of rings having almost minimal multiplicity (and maximal type) whose associated graded rings are not Cohen-Macaulay \cite{S4}.

Based on her results, Sally raised a conjecture predicting that Cohen-Macaulay local rings with almost minimal multiplicity have almost Cohen-Macaulay associated graded rings, i.e.,
 the depth is at least $d-1$, where $d={\rm dim}\,R$ \cite{S4}. This is known as {\em Sally's conjecture}. It was proved independently by Rossi-Valla \cite{RV1} and Wang \cite{W} about 13 years later.
However, it was not known whether
the associated graded ring is almost Cohen-Macaulay, provided the stretchedness of $\m$. In this direction, Rossi and Valla generalized the concept of stretchedness to the case of $\m$-primary ideals and proved effective criteria for the (almost) Cohen-Macaulayness
of associated graded rings  of stretched $\m$-primary ideals  \cite{RV3}. However, since the definition of stretched $\m$-primary ideals does not 
include  ideals of almost minimal multiplicity, these results cannot be viewed as a generalized version of Sally's conjecture. 

Recently, Polini and Xie extended the concepts of minimal and almost minimal multiplicity to ideals which are not necessarily primary to the maximal ideal by introducing the notion of minimal and almost minimal $j$-multiplicity \cite{PX1}. Under certain residual conditions, they proved that the associated graded rings of ideals having minimal $j$-multiplicity (almost minimal $j$-multiplicity, respectively) are always Cohen-Macaulay (almost Cohen-Macaulay, respectively).
 Their result on the almost Cohen-Macaulayness of associated graded rings extends for the first time Sally's conjecture to a class of ideals of arbitrary height, namely ideals having  almost minimal $j$-multiplicity.

In the present paper,  we generalize all the above results by proving a characterization and a sufficient condition, respectively, for the Cohen-Macaulayness and almost Cohen-Macaulayness of the associated graded rings of a new  class of ideals, dubbed $j$-stretched ideals. The concept of $j$-stretched ideals we introduce  includes all the above definitions given by  Rossi-Valla and Polini-Xie, i.e., stretched $\m$-primary ideals,
 ideals having  minimal $j$-multiplicity and ideals having almost minimal $j$-multiplicity.
In the $0$-dimensional case (that is, when both notions of stretchedness are defined) $j$-stretched ideals generalize stretched ideals in a non trivial way. Moreover, we exhibit several classes of examples over $1$-dimensional Cohen-Macaulay local rings that are $j$-stretched, they even have almost minimal multiplicity, but are not stretched
(see Examples \ref{nstr}, \ref{nstr2} and \ref{nstr3}).
 
 Our work, then, unifies the approaches of Rossi-Valla and Polini-Xie, and provides
  an extended  version of Sally's conjecture  for  this wider class of ideals.
 
Inspired by \cite{PX1}, our basic tools to study ideals of arbitrary height are general elements,  residual intersection theory (a generalization of linkage), and $j$-multiplicity theory (a higher dimensional version of the Hilbert-Samuel multiplicity). However, to deal with stretchedness and index of nilpotency,
one has to reduce to the case of finite length  by factoring out a sequence of  elements. The problem arises because  the length depends on the choice of the sequence of elements. To overcome this difficulty, we prove a crucial `Specialization Lemma' (Lemma \ref{specializ}). It states that if we factor out a sequence of general elements, we obtain a fixed length (which is also a lower bound for the length obtained by factoring out any special sequence with the same number of elements).

Lemma \ref{specializ} is employed in the proof of many results of the present paper. Moreover, it seems to have many applications in the process of generalizing classical results on multiplicities and Hilbert functions to arbitrary ideals and modules.

The structure of the paper is the following: in Section 2,  we define the concept of $j$-stretched ideals and show that it naturally extends  the notions of minimal and almost minimal $j$-multiplicity.
Thanks to the Specialization Lemma (Lemma \ref{specializ}), we prove that this concept also generalizes the classical definition of stretched $\m$-primary ideals (Theorem \ref{AN}, Corollary \ref{mprim}). We provide several classes of examples showing that in general, $j$-stretched $\m$-primary ideals need not be stretched (Examples \ref{nstr}, \ref{nstr2}, \ref{nstr3} and \ref{nstr4}). However, in contrast to these examples, we provide a sufficient condition for these two notions to coincide (Proposition \ref{equiv}). Finally, we answer  a question of  Sally that asks how does the stretchedness depend on minimal reductions (Corollary \ref{Sally}).
Section 3 is rather technical and collect tools needed to prove the main results of this paper. In this section we  describe  the behavior of general minimal reductions with respect to some specific properties (e.g., the nilpotency index), and the structure of $j$-stretched ideals. In Section 4,  we prove a clean characterization of the Cohen-Macaulayness of the associated graded rings of $j$-stretched ideals (Theorem \ref{CM}). We also prove a sufficient condition for the almost Cohen-Macaulayness of the associated graded ring of $j$-stretched ideals (Theorem \ref{2}). This latter result is a proof of Sally's conjecture for this class of ideals. Finally, we provide several applications of these theorems. Among the others, we recover the main results of \cite{PX1} and \cite{RV3}, and prove, under some additional assumptions, that the associated graded rings of ideals with almost almost minimal $j$-multiplicity are almost Cohen-Macaulay (Corollary \ref{almalm}).

\section{$j$-stretched ideals: basic properties}

In this section we  fix  notation and recall some basic facts that will be used throughout the paper. We then give the definition of $j$-stretched ideals and  show that this concept generalizes stretched $\m$-primary ideals,  ideals having minimal $j$-multiplicity and ideals having almost minimal $j$-multiplicity.

Let $(R, \m)$ be a Noetherian local ring with infinite residue field $k$ (after possibly replacing $R$ by $R(X)=R[X]_{\m R[X]}$, where $X$ is a variable over $R$, we can always assume the residue field of $R$ is infinite). Let $M$ be a  finite  $R$-module  and $I$ an arbitrary
$R$-ideal. Recall that $G={\rm
gr}_I(R):=\oplus_{j=0}^{\infty}I^j/I^{j+1}$ is the {\it associated graded ring} of $I$ and $T={\rm gr}_I(M):=\oplus_{j=0}^{\infty}I^jM/I^{j+1}M$ is the {\it associated graded module} of $I$ on $M$.
An $R$-ideal $J\subseteq I$ is said to be a {\em reduction} of $I$ on $M$ if there exists a non-negative integer $r$ such that $I^{r+1}M=JI^rM$. If further $J$ does not contain properly any other reduction of $I$ on $M$, we say that $J$ is a {\em minimal reduction} of $I$ on $M$. Since $|k|=\infty$, minimal reductions of $I$ on $M$ always exist. Furthermore, every minimal reduction of $I$ on $M$ can be generated minimally on $M$ by the same number of generators, dubbed the {\em analytic spread} $\ell(I,M)$ of $I$ on $M$.

Let $I =(a_1, \ldots,a_n)$ and write $x_i =\sum_{j=1}^n\lambda_{ij}a_j$ for $i=1,\dots,t$ and $(\lambda_{ij})\in R^{tn}$. The elements $x_1,\dots,x_t$ are said to form a {\em sequence of general elements} in $I$ (or, equivalently, $x_1,\dots,x_t$ are {\em general} in $I$) if there exists a dense open subset $U$ of $k^{tn}$ such that the image $(\overline{\lambda_{ij}}) \in U$. When $t=1$ one says that $x_1=x$ is a {\em general element} in $I$.
The relevance of this notion in our analysis comes from the following facts: (a) general elements always form a superficial sequence for $I$ on $M$ (\cite[Corollary~2.5]{X}); (b) if $t=\ell(I,M)$,  $x_1,\dots,x_t$ form a minimal reduction of $I$ on $M$ with reduction number  $r(I,M)$  (see for instance \cite[Corollary~2.2]{T2}); (c) one can compute the $j$-multiplicity using a general minimal reduction (\cite[Proposition~2.1]{PX1}).

Assume ${\rm dim}\,M=d$. One says that $I$ has {\em maximal analytic spread} on $M$ if $\ell(I, M)=d$.
In this case, write $J=(x_1, \ldots, x_d)$ for a general minimal reduction of $I$ on $M$, $J_{d-1}=(x_1,\ldots,x_{d-1})M:_M I^{\infty}$ and  $\overline{M}=M/J_{d-1}$. Then, $\overline{M}$ is  a 1-dimensional Cohen-Macaulay $R$-module and $I$ is an {\it ideal of definition} on $\overline{M}$, i.e.,  $\lambda(\overline{M}/I\overline{M})<\infty$, where $\lambda (N)$ denotes the length of $N$. Therefore, one can define the Hilbert function
of $I$ on $\overline{M}$ as follows
$$
HF_{I,\,\overline{M}}(j)=\lambda_{R}(I^j\overline{M}/I^{j+1}\overline{M}).
$$
The $j$-multiplicity (see for instance \cite[2.1]{PX1}) is
$$
j(I, M)=e(I, \overline{M})=\lambda(\overline{M}/x_d \overline{M})=\lambda(I\overline{M}/x_d I\overline{M})=\lambda(I\overline{M}/I^2\overline{M})+\lambda(I^2\overline{M}/x_dI\overline{M}).
$$
It is proved in \cite{PX2}  that the Hilbert function of $I$ on $\overline{M}$ does not depend on the general minimal reduction $J$.


We are ready to give the definition of   $j$-{\rm stretched} ideals on $M$.
\begin{Definition}\label{Def}
Let $M$ be a finite module of dimension $d$ over a Noetherian local ring and  $I$  an ideal with $\ell(I,M)=d$.
We say that $I$ is $j$-{\rm stretched} on $M$  if for a general minimal reduction $J=(x_1, \ldots, x_{d})$ of $I$ on $M$ one has
$$\lambda(I^2\overline{M}/x_dI\overline{M}+I^3\overline{M})\leq 1,$$
where $\overline{M}=M/J_{d-1}$ and $J_{d-1}=(x_1,\ldots,x_{d-1})M:_M I^{\infty}$.
\end{Definition}

Definition \ref{Def} relies on  general minimal reductions. However,  we will show later that under some reasonable (technical) assumptions, one can use any  reduction to prove the $j$-stretchedness of $I$ on $M$ (see Theorem \ref{AN}).

Assume $I$ is $j$-stretched on $M$. Set $t-1=\lambda(I^2\overline{M}/x_d\overline{M}\cap I^2\overline{M})$. It is easily seen that
$HF_{I,\,\overline{M}/x_d\overline{M}}(j)=1$ for $2 \leq j\leq t$, and $HF_{I,\,\overline{M}/x_d\overline{M}}(j)=0$ for $ j\geq t+1$. In other words, if $I$ is $j$-stretched on $M$, then $\overline{M}/x_d\overline{M}$ is a stretched Artinian module.

Clearly, $j$-stretched ideals include  ideals having minimal or almost minimal $j$-multiplicity, whose definitions are now recalled (see \cite{PX1}).

\begin{itemize}
\item I has {\em minimal $j$-multiplicity on $M$} if $\lambda(I^2\overline{M}/x_dI\overline{M})=0$;
\item I has {\em almost minimal $j$-multiplicity on $M$} if $\lambda(I^2\overline{M}/x_dI\overline{M})\leq 1$.
\end{itemize}
In particular, every $\m$-primary ideal having minimal or almost minimal multiplicity is $j$-stretched.

\begin{Example}\label{amm}$($see \cite{NU} and \cite[Example~4.10]{PX1}$)$
Let $S$ be a $3$-dimensional Cohen-Macaulay local ring and $x$,$y$, $z$ a system of parameters
for $S$. Set $R = S/(x^2 -yz)S$ and $I = (x,y)R$. Then $I$ is $j$-stretched on $R$.
\end{Example}

\demo This is clear because $I$ has minimal $j$-multiplicity ( \cite{PX1}).
\QED
\bigskip

The next result shows that $j$-stretchedness is preserved by faithfully flat extensions.
\begin{Lemma}\label{ff}
Let $M$ be a finite module of dimension $d$ over a Noetherian local ring $(R,\m)$ and $I$  an $R$-ideal. Assume $I$ is $j$-stretched on $M$ and $(S,\n)$ is a Noetherian local ring that is flat over $R$ with $\m S=\n$. Then $IS$ is $j$-stretched on $M\otimes_R S$.
\end{Lemma}

\demo
Set $M'=M\otimes_R S$. It is clear from  the dimension formula  that ${\rm dim}_R\,M={\rm dim}_S\,M'$. One also has $\ell(I,M)=\ell(IS, M')$ (see for instance \cite[Proposition~1.5]{AM}).
Now, set
$$L=I^2M/(x_dIM+I^3M+[((x_1,\dots,x_{d-1})M:_MI^{\infty}) \cap I^2M]).$$
 By assumption,  $\lambda_R(L)\leq 1$. Hence by the flatness and $\m S=\n$,
\begin{equation}\label{otimes}
\lambda_S(L\otimes_R S)=\lambda_R(L)\cdot \lambda(S/\m S)\leq \lambda(S/\m S)\leq 1.
\end{equation}
Again by flatness,
$$L\otimes_R S=I^2M'/(x_dIM'+I^3M'+[((x_1,\dots,x_{d-1})M^{\prime}:_{M'}I^{\infty}) \cap I^2M']).$$
This fact, together with (\ref{otimes}), shows that $IS$ is $j$-stretched on $M'=M\otimes_R S$.
\QED
\bigskip

As a consequence of Lemma \ref{ff}, we obtain immediately that the property of being $j$-stretched is preserved by passing to the completion of $R$, or enlarging the residue field, that is, replacing $R$ by $R(X)=R[X]_{\m R[X]}$, where $X$ is a variable over $R$.

Next, we want to show that the $j$-stretchedness generalizes the notion of stretched $\m$-primary ideals given by Rossi and Valla in \cite{RV3}.
  Recall that $I$ is a  {\em stretched} $\m$-primary ideal if $I$ is $\m$-primary and there exists a minimal reduction $J$ of $I$ such that
\begin{itemize}
\item[(a)]  $J\cap I^2=JI$;
\item[(b)]  $HF_{I/J}(2)\leq 1$.
\end{itemize}

The notion of stretched $\m$-primary ideals extends the classical definition of stretched maximal ideals given by Sally in \cite{S3}.
If $R$ is Cohen-Macaulay, it also includes ideals having minimal multiplicity. However,
we found several classes of $\m$-primary ideals with almost minimal multiplicity that are not stretched, even in $1$-dimensional Cohen-Macaulay local rings - see for instance Examples \ref{nstr}, \ref{nstr2} and \ref{nstr3}.
Hence, stretched $\m$-primary ideals do not generalize ideals with almost minimal multiplicity.

On the other hand, the notion of $j$-stretchedness does not present this pathology, because it generalizes naturally the concepts of minimal and almost minimal $j$-multiplicity (which, in turn, generalize ideals having minimal multiplicity and ideals having almost minimal multiplicity, respectively).
\smallskip

Notice that stretched ideals require the existence of {\em one} minimal reduction $J$ of $I$ satisfying properties (a) and (b) above, whereas the definition of $j$-stretched ideals imposes a condition on every {\em general} minimal reduction. Hence, {\em a priori}, it is not clear how these two notions relate. To answer this natural question we prove a crucial specialization lemma. Before stating it, we recall
 the notion of specialization of modules, introduced by Nhi and Trung \cite{NT}.

Let $S=R[\underline{z}]$, where  $ \underline{z}=z_1, \ldots,  z_{t}$ are variables  over the Noetherian local ring $(R,\m)$. Let $M^{\prime}$ be a finite $S$-module.
 Let $\phi:
S^f \rightarrow S^g \rightarrow  0$ be a finite free presentation of
$M^{\prime}$ and let $A=(a_{ij}[\underline{z}])$ be a matrix representation of $\phi$.
For any vector $\underline{\alpha}=(\alpha_1, \ldots, \alpha_t)\in R^t$, let $A_{\underline{\alpha}}:=(a_{ij}[\underline{\alpha}])$ and $\phi_{\underline{\alpha}}: R^f\rightarrow R^g\rightarrow  0$ be the corresponding map defined by $A_{\underline{\alpha}}$. One says that $\phi_{\underline{\alpha}}$ is a {\it specialization} of $\phi$. A {\it specialization} of $M^{\prime}$ is defined to be $M^{\prime}_{\underline{\alpha}}:={\rm Coker}(\phi_{\underline{\alpha}})$. By \cite{NT}, $M^{\prime}_{\underline{\alpha}}$ depends only on $\underline{\alpha}$ (it does not depend on the choice of $\phi$ and $A$). The vector $\underline{\alpha}\in R^t$ is said to be {\it general} (equivalently, the specialization $M^{\prime}_{\underline{\alpha}}$ is {\it general}) if the image $\overline{\underline{\alpha}}=(\overline{\alpha}_1, \ldots, \overline{\alpha}_t)\in U$ where $U$ is some dense open subset of $k^t$, where $k=R/\m$.

\begin{Lemma}\label{specializ}$($Specialization Lemma$)$
Let $R, S$ be defined as above. Assume $k$ is infinite.
 Let $M$ be a finite $R$-module and $M^{\prime}=M\otimes_R S$.
  Let $N^{\prime} \subseteq M^{\prime}$ be a submodule
such that
$\lambda_{S_{\m S}}(M^{\prime}_{\m S}/N^{\prime}_{\m S})=\delta$.  Then
\begin{itemize}
\item[(a)]  for a general vector $\underline{\alpha} \in R^t$,  the length $\lambda_R(M/N^{\prime}_{\underline{\alpha}})=\delta$.
\item[(b)]   Assume $R$ is equicharacteristic and fix any vector $\underline{\alpha_0}\in R^t$, Then for a general vector $\underline{\alpha}  \in R^t$,  the length $\delta=\lambda_R(M/N^{\prime}_{\underline{\alpha}})\leq \lambda_R(M/N^{\prime}_{\underline{\alpha_0}})$.
\end{itemize}
\end{Lemma}

\demo
We may pass to the $\m$-adic completion of $R$ to assume that $R$ contains
a regular local ring and therefore $M$ is a finite module over a polynomial ring over a regular local ring.
We use induction on $\delta$ to prove part (a). Clearly, this statement holds if $\delta=0$. We then assume that $\delta>0$ and assertion (a) holds for $\delta-1$. Let $$M^{\prime}=M_0^{\prime}  \supset  M_1^{\prime}\supset \ldots \supset M_{\delta}^{\prime}=N^{\prime}$$
be a filtration such that $(M_{l-1}^{\prime}/M_{l}^{\prime})_{\m S}\cong k(\underline{z})$
for $1\leq l\leq \delta$. By induction hypothesis,  for a
general  $\underline{\alpha}\in R^t$ we have
 $\lambda_{R}(M/(M^{\prime}_{\delta-1})_{\underline{\alpha}})= \delta-1$.
Applying an argument similar to the one of \cite{NT}, for a general  $\underline{\alpha}\in R^t$ we obtain $$(M^{\prime}_{\delta-1})_{\underline{\alpha}}/(M^{\prime}_{\delta})_{\underline{\alpha}}\cong (M_{\delta-1}^{\prime}/M_{\delta}^{\prime})_{\underline{\alpha}}\cong k(\underline{z})_{\underline{\alpha}}\cong k.$$
Hence, for a general $\underline{\alpha}\in R^t$, we obtain
$$\lambda_{R}(M/N^{\prime}_{\underline{\alpha}})=\lambda_{R}(M/(M^{\prime}_{\delta-1})_{\underline{\alpha}})+
\lambda_{R}((M^{\prime}_{\delta-1})_{\underline{\alpha}}/(M^{\prime}_{\delta})_{\underline{\alpha}})=\delta.$$

To prove part (b), since $R$ is equicharacteristic,  after possibly passing to the $\m$-adic completion of $R$, we may assume that $R$ contains its residue field $k$.
For every $\underline{\alpha}$ in $k^t$, we have the isomorphisms
$$M^{\prime}/N^{\prime}\otimes_{S} S/(\m, \underline{z}-\underline{\alpha}) \cong M^{\prime}/N^{\prime} \otimes_{S}
S/\m S \otimes_{k[\underline{z}]} k[\underline{z}]/(\underline{z}-\underline{\alpha})
$$
$$\cong  M^{\prime}/N^{\prime} \otimes_{k[\underline{z}]} k[\underline{z}]/(\underline{z}-\underline{\alpha})\cong
M/N^{\prime}_{\underline{\alpha}},$$
where the second isomorphism follows because $M^{\prime}/N^{\prime}$ is a $S/\m S$ module.
Therefore, we have
$$\mu_{S_{(\m, \underline{z}-\underline{\alpha})S}}(M^{\prime}/N^{\prime} \otimes_{S} S_{(\m, \underline{z}-\underline{\alpha})S})=\lambda_{S_{(\m, \underline{z}-\underline{\alpha})S}}(M^{\prime}/N^{\prime} \otimes_{S} S/(\m,\underline{z}-\underline{\alpha})S)=\lambda_R(M/N^{\prime}_{\underline{\alpha}}),$$
because $k$ is the residue field of $S_{(\m, \underline{z}-\underline{\alpha})S}$ and is contained in $S_{(\m, \underline{z}-\underline{\alpha})S}$.

Now, set $q=\lambda_R(M/N^{\prime}_{\underline{\alpha_0}})$. Then,
$$U=\{ \underline{\alpha} \in k^t \,| \, \mu_{S_{(\m, \underline{z}-\underline{\alpha})S}} (M^{\prime}/N^{\prime} \otimes_{S} S_{(\m, Z- \Lambda)S})\leq q\}= k^t \setminus V({\rm Fitt}_q(M^{\prime}/N^{\prime})) $$
is a Zariski open subset of $k^t$. Since $\underline{\alpha_0} \in U$, we also know it is non-empty and hence it is dense.
Then, for any $\underline{\alpha} \in U$, we have
$$\lambda_R(M/N^{\prime}_{\underline{\alpha}} )= \mu_{S_{(\m, \underline{z}-\underline{\alpha})S}}(M^{\prime}/N^{\prime} \otimes_{S} S_{(\m,\underline{z}-\underline{\alpha})S}) \leq  q=\lambda_R(M/N^{\prime}_{\underline{\alpha_0}}).$$

\QED
\bigskip

Lemma \ref{specializ} greatly enhances our ability to deal with  arbitrary ideals.
In the case of non $\m$-primary ideals, to get a finite length one has to factor out a sequence of  elements. But  the  length usually depends on the choice of the elements in the sequence. Lemma \ref{specializ} shows that this issue is fixed if one restricts to sequences of general elements.

We now recall some definitions from the theory of residual intersections (see, for instance, \cite{U} or \cite{PX1}). Let $M$ be a finite faithful module over a Noetherian local ring $R$, and $I$ be an ideal and $x_1,\ldots,x_t$ be elements of $I$. 
Write
$H=(x_1, \ldots, x_t)M:_M I$. If $IM_\p=(x_1, \ldots, x_{t})M_\p$
for every $\p\in {\rm Spec}(R)$ with ${\rm dim}\,R_\p\leq t-1$, then
$H$ is said to be a {\it  $t$-residual intersection } of $I$ on $M$.
Now let $H$ be a $t$-residual intersection of $I$ on $M$. If in
addition $IM_\p=(x_1, \ldots, x_{t})M_\p$ for every  $\p\in {\rm
Supp}_R(M/IM)$ with ${\rm dim}\,R_\p\leq t$, then $H$ is said to be
a {\it geometric $t$-residual intersection } of $I$ on $M$. If $M$
is not faithful, then we say that $I$  satisfies the {\it condition
$G_s$} on $M$ if $I\overline{R}$  satisfies the condition $G_s$ on
$M$, where $\overline{R}=R/{\rm Ann}\,M$. We then say $H$ is a {\it
$t$-residual intersection } or {\it geometric $t$-residual
intersection } of $I$ on $M$ if $H$ is a $t$-residual intersection
or geometric $t$-residual intersection of $I\overline{R}$ on $M$
respectively.

If  $M$ is a finite module over a catenary and equidimensional Noetherian local ring  and $I$ is an ideal which satisfies $G_s$ on $M$,  then for general elements $x_1, \ldots, x_s$ in $I$ and $0\leq i< s$, $H_i=(x_1, \ldots, x_i)M:_M I$ is a geometric $i$-residual intersection of $I$ on $M$ and $H_s=(x_1, \ldots, x_s)M:_M I$ is a $s$-residual intersection of $I$ on $M$ \cite[Lemma~3.1]{PX1}.
Finally, let  $M$ be   Cohen-Macaulay, then $I$ is said to have the {\it Artin-Nagata property} ${AN^-_t}$ on $M$ if for every $0\leq i\leq t$ and every geometric $i$-residual intersection $H$ of $I$ on $M, \,$  the module $M/H$ is  Cohen-Macaulay.

In the next result we employ the Specialization Lemma to prove that  general minimal reductions achieve the minimal colength.
\begin{Proposition}\label{generic1}
Let $M$ be a finite module of dimension $d$ over a Noetherian local ring $R$ with infinite residue field. Let $I$ be an ideal with $\ell(I,M)=d$. Assume $I$ satisfies $G_d$ on $M$. Let $H$ and $J$ be a minimal and a general minimal reduction of $I$ on $M$, respectively.  Let  $j\geq 1$ be a fixed integer. Then the lengths $\lambda(I^jM/J^jM)$ and $\lambda(I^jM/JI^{j-1}M+I^{j+1}M)$ do not depend on $J$. Furthermore if $R$ is  equicharacteristic  then
\begin{itemize}
\item[(a)] $\lambda(I^jM/J^jM)\leq \lambda(I^jM/H^jM)$.
\item[(b)] $\lambda(I^jM/JI^{j-1}M+I^{j+1}M)\leq \lambda(I^jM/HI^{j-1}M+I^{j+1}M)$.
\end{itemize}
\end{Proposition}

\demo Let $\m$ be the maximal ideal of $R$  and write $I=(a_1, \ldots, a_n)R$.
We first prove (a). Assertion (b) can be done similarly. Take $d \times n$ variables, say $\underline{z}=(z_{ij})$, and set $S=R[\underline{z}]$, $M^{\prime}=M\otimes_{R}S$, $J^{\prime}=(x_1^{\prime}, \ldots x_{d}^{\prime})S$, where $x_i^{\prime}=\sum_{j=1}^n z_{ij}a_j$, $1\leq i \leq d$.  Let  $\underline{\alpha_0}\in R^{dn}$ be the vector such that $J^{\prime}_{\underline{\alpha_0}}=H$.
Since $I$ satisfies $G_d$ on $M$, we have $\lambda_{S_{\m S}}(IM^{\prime}_{\m S}/J^{\prime}M^{\prime}_{\m S})<\infty$.
By Lemma \ref{specializ}, for a general element $\underline{\alpha}\in R^{dn}$,
$$\lambda(I^jM/J^jM)=\lambda(I^jM/(J^{\prime})^jM^{\prime}_{\underline{\alpha}})=
\lambda(I^jM^{\prime}/(J^{\prime})^jM^{\prime}).$$
Furthermore, if $R$ is  equicharacteristic,
$$\lambda(I^jM/J^jM)\leq
\lambda(I^jM/(J^{\prime})^jM^{\prime}_{\underline{\alpha_0}})=\lambda(I^jM/H^jM).$$
\QED
\bigskip

A first consequence of Proposition \ref{generic1} is the following result, describing the behaviors of intersections with respect to a minimal reduction and a general minimal reduction.

\begin{Proposition}\label{inters}
Let $M$ be a Cohen-Macaulay module of dimension $d$ over an equicharacteristic Noetherian local ring and $I$ an ideal with  $\lambda(M/IM)<\infty$.
If $H$ is a minimal reduction of $I$ on $M$ and $J$ is a general minimal reduction of $I$ on $M$, then
$$\lambda((JM \cap I^2M   )/JIM)\leq \lambda((HM\cap I^2M )/HIM).$$
In particular, if $HM\cap I^2M =HIM$ then $JM \cap I^2M =JIM$.
\end{Proposition}

\demo Since  for any ideal $K\subseteq I$,
$$\lambda((KM \cap I^2M )/KIM)=\lambda(I^2M/KIM)-\lambda(I^2M/KM\cap I^2M)$$
and
$$\lambda(I^2M/KM\cap I^2M )=\lambda(M/KM)-\lambda(M/(I^2M+KM))$$
 we have
\begin{equation}\label{intersez}
\lambda((KM \cap I^2M )/KIM)=\lambda(I^2M/KIM)-\lambda(M/KM)+\lambda(M/(I^2M+KM)).
\end{equation}
Now, observe that \begin{itemize}
\item $\lambda(I^2M/JIM)=\lambda(I^2M/HIM)$ (see, for instance, \cite[Corollary~2.1]{RV});
\item $\lambda(M/JM)=e(M)=\lambda(M/HM)$, because $J$ and $H$ are minimal reductions of $I$ on $M$;
\item $\lambda(M/(I^2M+JM))\leq \lambda(M/(I^2M+HM))$, by Lemma \ref{specializ}.
\end{itemize}
These three facts, together with (\ref{intersez}), show that
$$\lambda((JM\cap I^2M)/JIM)\leq \lambda((HM\cap I^2M)/HIM).$$
\QED
\bigskip

We are now ready to prove the main result of this section. It shows that, in the $\m$-primary case or in presence of residual assumptions, $j$-stretchedness can be proved by checking {\em any} minimal reduction (instead of of checking every general minimal reduction).
\begin{Theorem}\label{AN}
Let $M$ be a Cohen-Macaulay module of dimension $d$ over an equicharacteristic Noetherian local ring $R$ and $I$ an $R$-ideal with $\ell(I,M)=d$. Let $H=(y_1,\dots,y_d)$ be a minimal reduction of $I$ on $M$. Set $H_{d-1}=(y_1,\ldots,y_{d-1})$ and assume
$$\lambda(I^2M/[y_dIM + I^3M + (H_{d-1}M:_M I^{\infty}) \cap I^2M])\leq 1.$$
If one of the two following conditions holds:
\begin{itemize}
\item[(i)] $\lambda(M/IM)<\infty$\, and\, $HM \cap I^2M=HIM;$
\item[(ii)] ${\rm depth}\,(M/IM)\geq 1$,  $I$ satisfies $G_d$ and $AN^-_{d-2}$ on $M$ and $H_{d-1}M:_M I$ is a geometric $d-1$-residual intersection of $I$ on $M;$
\end{itemize}
then $I$ is $j$-stretched on $M$.

\end{Theorem}

From the computational point of view, Theorem \ref{AN} is useful to produce examples of $j$-stretched ideals.
Indeed, it says that, if we can find one minimal reduction of $I$ (on $M$) satisfying a certain inequality, we have immediately that $I$ is $j$-stretched (on $M$).
\medskip

\demo
Since  $H_{d-1}M:_M I$ is a geometric $d-1$-residual intersection of $I$ on $M$ and $I$ satisfies $AN^-_{d-2}$ on $M$, we have
$$H_{d-1}M:_MI^{\infty}=H_{d-1}M:_MI.$$
We want to show that in our setting, either (i) or (ii) imply $(H_{d-1}M:_M I^{\infty}) \cap I^2M = H_{d-1}IM$, yielding immediately that
$$ \lambda(I^2M/(HIM + I^3M))\leq 1.$$
First assume (i) holds. Since $M$ is Cohen-Macaulay, $I$ contains a non zero divisor on $M/H_{d-1}M$. Hence
$H_{d-1}M:_M I^{\infty} = H_{d-1}M$, yielding that
$$(H_{d-1}M:_M I^{\infty}) \cap I^2M = H_{d-1}M \cap I^2M = H_{d-1}M \cap HM \cap I^2M = H_{d-1}M \cap HIM.$$
Now,
$$
\begin{array}{ll}
H_{d-1}M \cap HIM & = H_{d-1}M \cap (H_{d-1}IM + y_dIM)\\
                                   & = H_{d-1}IM + (H_{d-1}M \cap y_dIM)\\
                                   & = H_{d-1}IM + y_d(H_{d-1}M :_{IM} y_d)\\
                                   & = H_{d-1}IM  + y_dH_{d-1}M\\
                                   & = H_{d-1}IM.
\end{array}
$$
Secondly, assume (ii) holds. An argument similar to \cite[Lemma~3.2.(f)]{PX1} gives the desired claim.
Hence, we have showed that in either case one has  $$\lambda(I^2M/HIM + I^3M)\leq 1.$$ Now, Lemma \ref{specializ} implies that
$\lambda(I^2M/JIM + I^3M)\leq 1$ for a general minimal reduction $J$ of $I$ on $M$. In turn, this yields that $I$ is $j$-stretched on $M$.
\QED
\bigskip

Taking $M=R$ in  Theorem \ref{AN} one immediately obtains that  every stretched $\m$-primary ideal is $j$-stretched.
\begin{Corollary}\label{mprim}
Let $(R,\m)$ be an equicharacteristic Cohen-Macaulay local ring and $I$ an $\m$-primary ideal. If $I$ is stretched, then $I$ is $j$-stretched on $R$.
\end{Corollary}

One may wonder if, in the $\m$-primary case, $j$-stretchedness coincides with stretchedness.
This is definitely not the case.
Indeed, we now present {\em classes} of ideals that are $j$-stretched (even have almost minimal multiplicity!) but are not stretched. This shows that $j$-stretched $\m$-primary ideals widely generalize stretched ideals already in the $1$-dimensional case.
\begin{Example}\label{nstr}
Fix $n\geq 3$. Let $A=k[\![t^n,t^{n+1},\ldots,t^{2n-1}]\!]$, $\m=(t^n,t^{n+1},\ldots,t^{2n-1})A$ and $I=(t^n,t^{n+1},\ldots,t^{2n-2})A$. Then $I$ is an $\m$-primary ideal that is $j$-stretched (has almost minimal multiplicity) but is not stretched with respect to {\rm any} minimal reduction.
\end{Example}

\demo
Notice that $I$ is $\m$-primary because it is a non zero ideal in a $1$-dimensional local domain. Next, we want to prove that $I$ has almost minimal multiplicity. To do it, we first observe that $H=(t^n)A$ is a minimal reduction of $I$. Indeed, $I^2=(t^{2n},t^{2n+1},\ldots,t^{3n-1})A$ and  it easy to check that $HI^2=(t^{3n}, t^{3n+1}, \ldots, t^{4n-1})A=I^3$. We now show that $\lambda(I^2/HI+I^3)=1$. First of all, notice that $I^2=H\m$. In particular, $I^2\subseteq H$, showing that $I^3\subseteq HI$. Hence, we only have to show that  $\lambda(I^2/HI)=1$.

Since $HI=(t^{2n},t^{2n+1},\ldots,t^{3n-2})A$, and $n\geq 3$, we have that $t^{3n-1} \in I^2 \setminus HI$, showing that
\begin{equation}\label{one} \lambda(I^2/HI)\geq 1.\end{equation}
The other generators of $I^2$, that is, $t^{2n},t^{2n+1},\ldots,t^{3n-2}$, are all in $HI$. Hence, to show that equality holds in (\ref{one}), we only need to check that $t^{3n-1}\m \subseteq HI$. However, it is easily seen that $t^{3n+k} \in HI$ for every $k\geq 0$, in particular, this proves that $t^{3n-1}\m \subseteq HI$.
Therefore, we have proved that
$$\lambda(I^2/HI)=1.$$
Now, by \cite[Corollary~2.1]{RV},
we obtain $\lambda(I^2/JI)=1$ for every minimal reduction $J$ of $I$. Therefore, $I$ has almost minimal multiplicity.

Next, we want to show that $I$ is not stretched with respect to any minimal reduction. Let $J$ be any principal reduction of $I$, we need to show that $J\cap I^2 \neq JI$. For a reduction $L$ of $I$, set
$$s_L(I)={\rm min}\{n \in \mathbb N \,| \, I^{n+1} \subseteq L\}$$
($s_L(I)$ is called the {\em index of nilpotency} of $I$ in $L$).
Since $I\neq H$ and $I^2=H\m \subseteq H$, it follows that
$$s_H(I)= 1.$$
Also, by the above, we have $I^2 \neq HI$, showing that $r_H(I)\geq 2$. However, we mentioned already that $I^3=HI^2$, which in turn proves $r_H(I)=2$.
Since we are in the $1$-dimensional case, $r_J(I)=2$ for every minimal reduction $J$ of $I$. Also, the trivial inequality $s_J(I)\leq r_J(I)$, shows that
$$s_J(I)\leq 2 \mbox{ for every minimal reduction } J \mbox{ of }I.$$
We claim that $s_J(I)$ must be exactly $1$ for every minimal reduction $J$ of $I$. Indeed, assume by contradiction that $s_{J'}(I)=2$ for some minimal reduction $J'$ of $I$. Then, for this minimal reduction $J'$ of $I$ we would have $s_{J'}(I)=r_{J'}(I)$. But then, by \cite[Remark~4.6.(2)]{HKU}, we would have
$$s_J(I)=r_J(I) \mbox{ for every minimal reduction } J \mbox{ of }I.$$
This is a contradiction, because we proved above that $s_H(I)\neq r_H(I)$. Therefore, we must have $s_J(I)=1$ for every minimal reduction $J$ of $I$.
Hence, we obtain that $I^2 \subseteq J$ for every minimal reduction $J$ of $I$.
Then, we have
$$\lambda(J\cap I^2/JI)=\lambda(I^2/JI)=1,$$
proving that $J\cap I^2 \neq JI$ for every principal reduction $J$ and showing that $I$ is not stretched with respect to {\em any} minimal reduction.

Finally, we need to show that $I$ is $j$-stretched on $R$. This follows either by the fact that $I$ has almost minimal $j$-multiplicity (because $I$ has almost minimal multiplicity), or, otherwise, by the following equality (obtained because $I^3\subseteq JI$)
$$\lambda(I^2/JI+I^3)= \lambda(I^2/JI)=1,$$
for any principal reduction $J$ of $I$.
\QED
\bigskip

If the maximal ideal of the semigroup ring has only three generators, we obtain more classes of examples.
\begin{Example}\label{nstr2}
Fix $n\geq 3$. Let $A=k[\![t^n,t^{n+1},t^{n+2}]\!]$, $\m=(t^n,t^{n+1},t^{n+2})A$ and $I=(t^n,t^{n+1})A$. Then $I$ is an $\m$-primary ideal
that is $j$-stretched, has almost minimal multiplicity, but is not stretched with respect to {\rm any} minimal reduction.
\end{Example}

\demo
Consider the minimal reduction $H=(t^n)A$ of $I$. Similar to the above, one can show that \begin{itemize}
\item $I^2=H\m \subseteq H$;
\item $r_H(I)=s_H(I)=2$;
\item $\lambda(I^2/HI)=1$.
\end{itemize}
As in the proof of Example \ref{nstr}, these three facts show that $r_J(I)=s_J(I)=2$ and $\lambda(I^2/JI)=1$, for any minimal reduction $J$ of $I$. Therefore, $I$ is $j$-stretched, has almost minimal multiplicity, but is not stretched.
\QED
\bigskip

Similar examples can be produced using arithmetic progressions. The proofs are similar to the previous ones.
\begin{Example}\label{nstr3}
Fix $n\geq 3$ and $a\geq 1$ with $2a=n$ or $3a=n$. Let $A=k[\![t^n,t^{n+a},t^{n+2a}]\!]$, $\m=(t^n,t^{n+a},t^{n+2a})A$ and $I=(t^n,t^{n+a})A$. Then $I$ is an $\m$-primary ideal
that is $j$-stretched, has almost minimal multiplicity, but is not stretched with respect to {\rm any} minimal reduction.
\end{Example}

As a consequence of either Example \ref{nstr}, \ref{nstr2} or \ref{nstr3}, one obtains that in the ring $A=k[\![t^3,t^{4},t^{5}]\!]$, the ideal $I=(t^3,t^4)$ is $j$-stretched, has almost minimal multiplicity, but is not stretched. These examples show that already in the $1$-dimensional case, there are many examples of $j$-stretched $\m$-primary ideals that are not stretched.

Next, we exhibit the example of an ideal that is $\m$-primary and $j$-stretched, but is not stretched and does not have almost minimal multiplicity.
\begin{Example}\label{nstr4}
Let $A=k[\![t^7,t^9,t^{10}]\!]$, $\m=(t^7,t^9,t^{10})A$ and $I=(t^7,t^9)A$. Then, $I$ is an $\m$-primary ideal that is $j$-stretched, is not stretched, and does not have almost minimal multiplicity.
\end{Example}

\demo Consider the minimal reduction $H=(t^7)A$ of $I$. Then, $$\lambda(I^2/HI+I^3)=1,$$
because $I^2=(t^{14},t^{16},t^{18})A$, $t^{14}$ and $t^{16}$ are in $HI$, $t^{18} \in I^2 \setminus HI+I^3$ and $t^{18}\m \subseteq HI+I^3$.
Similarly
 $$\lambda(I^3/HI^2+I^4)=1,$$
because $I^3=(t^{21},t^{23},t^{25},t^{27})A$, $t^{21}$, $t^{23}$ and $t^{25}$ are in $HI^2$, $t^{27}\in I^3 \setminus HI^2+I^4$ and $t^{27}\m \subseteq HI^2+I^4$.
Again by the same argument,
$$\lambda(I^4/HI^3+I^5)=1,$$
in this case, the key element is $t^{36}$.
One can also check that $$\lambda(I^k/HI^{k-1}+I^{k+1})=0, \mbox{ for every } k\geq 5.$$
Since $\lambda(I^2/HI+I^3)=1$, by Proposition 2.5, $I$ is $j$-stretched. Clearly $I$ does not have almost minimal $J$-multiplicity.
Finally, one can see that $I$ is not stretched, because it fails the intersection property.
\QED
\bigskip

Similar to Example \ref{nstr4}, one could prove that if $A=k[\![t^5,t^7,t^{8}]\!]$, $\m=(t^5,t^7,t^{8})A$ and $I=(t^5,t^7)A$. Then, $I$ is an $\m$-primary ideal that is $j$-stretched, is not stretched, and does not have almost minimal multiplicity. In this case, the $j$-multiplicity of $I$ is one less than the $j$-multiplicity of the ideal in Example \ref{nstr4}. 

In contrast to the previous examples, we are able to provide situations where $j$-stretchedness coincides with stretchedness. The next result, for instance, shows that, for $\m$-primary ideals, stretchedness and $j$-stretchedness are the same when the associated graded ring has some good properties.
\begin{Proposition}\label{equiv}
Let $R$ be a Cohen-Macaulay local ring and $I$ an $\m$-primary ideal. Assume $I^2\cap H=HI$ for a minimal reduction $H$ of $I$ $($e.g., if ${\rm gr}_I(R)$ is Cohen-Macaulay$)$.
Then, $I$ is stretched if and only if $I$ is $j$-stretched.
\end{Proposition}

\demo  Since one implication has been proved in Corollary \ref{mprim}, we only need to show that if $I$ is $j$-stretched, then $I$ is stretched. Let $J=(x_1,\dots,x_d)$ be a general minimal reduction of $I$ and $J_{d-1}=(x_1,\dots,x_{d-1})$. By $j$-stretchedness we have
 $$\lambda(\overline{I}^2/(x_d\overline{I}+\overline{I}^3)\leq 1,$$
 where $\overline{R}=R/(J_{d-1}:I^{\infty})=R/J_{d-1}$.
By assumption and Proposition \ref{inters}, one obtains that $J\cap I^2=JI$. Hence, we get
$$\overline{I}^2/(x_d\overline{I}+\overline{I}^3)\simeq I^2/(JI+I^3+J_{d-1}\cap I^2)= I^2/(JI+I^3)=I^2/((J\cap I^2)+I^3).$$
Therefore, for a general minimal reduction $J$ of $I$, we have $HF_{I/J}(2)\leq 1.$ This fact, together with $J\cap I^2=JI,$ shows the stretchedness of $I$.
\QED
\bigskip

We finish this section with an application of the above results to answer a question of Sally. If $I$ is $\m$-primary and $I\neq \m$, there are classical examples (provided by Sally or Rossi-Valla) showing that $I$ can be stretched with respect to a minimal reduction $J_1$ but not stretched with respect to a different minimal reduction $J_2$. Hence, it is well known that the stretchedness property depends upon the minimal reduction. Sally in \cite{S3} raised the following question: {\em to what extent does the concept of `stretchedness' depend upon the choice of the minimal reduction?}

Thanks to the Specialization Lemma, we are now able to answer this question.

\begin{Corollary}\label{Sally}
Let $(R,\m)$ be an equicharacteristic Cohen-Macaulay local ring and $I$ an $\m$-primary ideal. If $I$ is stretched with respect to a minimal reduction $H$, then $I$ is stretched with respect to any general minimal reduction.
\end{Corollary}

\demo Let $J$ be a general minimal reduction of $I$. By Proposition \ref{inters}, the `intersection property' $J\cap I^2=JI$ follows at once from $H\cap I^2=HI$.
Then, we only need to show that $\lambda(I^2/(J\cap I^2) + I^3)\leq 1$. By Proposition \ref{generic1} we have $\lambda(I^2/JI+I^3)\leq \lambda(I^2/HI+I^3)$. Hence, we obtain the following chain of inequalities
$$\lambda(I^2/(J\cap I^2) + I^3)=\lambda(I^2/JI+I^3)\leq \lambda(I^2/HI+I^3)=\lambda(I^2/(H\cap I^2)+I^3)\leq 1.$$
\QED
\bigskip

Corollary \ref{Sally} shows that, in an equicharacteristic Cohen-Macaulay local ring, the property of being a stretched ideal does not depend on the choice of the general minimal  reduction.
In some sense, this answer is the best that one could hope. Indeed, stretchedness depends heavily upon the choice of the minimal reduction, however Corollary \ref{Sally} states that `most' reductions exhibit the same behavior, since the property of being stretched holds for minimal reductions from a dense  open Zariski subset.
\smallskip

\section{Structure of $j$-stretched ideals}

In this section we prove  some technical results on numerical invariants of $j$-stretched ideals  that will be used to prove the main results of  next section.
\medskip

We  define the {\it index of nilpotency of $I$ on $M$} with respect to a reduction $J$  to be
$$
s_J(I, M)={\rm min}\{j\, |\, I^{j+1}M\subseteq JM\}.
$$

The following two results (Lemma \ref{nilpotency1} and Proposition \ref{nilpotency2}) generalize (to arbitrary ideals),
 Propositions \cite[5.3.2]{F} and \cite[5.3.3]{F} from Fouli's thesis, where she showed that for $\m$-primary ideals
  the index of nilpotency does not depend on the general minimal reduction,
  and general minimal reductions always achieve the largest possible index of nilpotency. 

The following lemma is the module-theoretic version  of   \cite[Proposition~5.3.2]{F}.

\begin{Lemma}$($\cite[5.3.2]{F}$)$ \label{nilpotency1}
Let $M$ be a finite module over a Noetherian local ring $(R, \m)$ of infinite residue field. Let $I$ be an $R$-ideal with $\ell(I, M)=\ell$.  Write $I=(a_1, \ldots, a_n)$, $S=R[\underline{z}]_{\m R[\underline{z}]}$, where $\underline{z}=(z_{ij})$ are $\ell \times n$ variables, $M^{\prime}=M\otimes_{R}S$, and $J^{\prime}=(x_1^{\prime}, \ldots, x_{\ell}^{\prime})S$, where $x_i^{\prime}=\sum_{j=1}^n z_{ij}a_j$, $1\leq i \leq \ell$. If  $J$ is a general minimal reduction of $I$ on $M$, then $s_J(I, M)\leq s_{J^{\prime}}(IS, M^{\prime})$.
\end{Lemma}

The next result generalizes \cite[Proposition~5.3.3]{F} to non $\m$-primary ideals. We need it in the module-theoretic version. Once again, the Specialization Lemma is a crucial ingredient in the proof.

\begin{Proposition}\label{nilpotency2}
Let $M$ be a Cohen-Macaulay module and $I$ an ideal over an equicharacteristic Noetherian local ring  with infinite residue field.
Let $H$ and $J$ be a minimal  and  a general minimal reduction of $I$ on $M$ respectively. One has the following statements:
 \begin{itemize}
\item[(a)] Assume $\lambda(M/IM)<\infty$. Then, $s_H(I, M)\leq s_J(I, M)$. In particular $s_J(I, M)$ does not depend on  $J$.
\item[(b)] Assume $\ell(I, M)={\rm dim}\,M=d$ and $I$ satisfies $G_d$ and $AN^-_{d-2}$ on $M$. Then $s_J(I, M)$ does not depend on  $J$.
    Furthermore, write $H=(y_1,\ldots,y_d)$ and assume $(y_1, \ldots, y_{d-2})M:_M I$ is a geometric $d-2$-residual intersection of $I$ on $M$ and $j(I,M)=e(I, \overline{M})$, where $\overline{M}=M/H_{d-1}$ and $H_{d-1}=(y_1, \ldots, y_{d-1})M:_M I^{\infty}$.
If\, ${\rm depth}\,(M/IM)\geq 1$ and $J\cap I^2M\subseteq JIM$, then $s_H(I, M)\leq s_J(I, M)$.
\end{itemize}
\end{Proposition}

\demo
We follow an argument similar to \cite[5.3.3]{F}. To prove part (a), consider the following two exact sequences:
\begin{eqnarray*}
& 0\rightarrow HM+I^{s+1}M/HM \rightarrow M/HM \rightarrow M/HM+I^{s+1}M \rightarrow 0,\\
& 0\rightarrow JM+I^{s+1}M/JM \rightarrow M/JM \rightarrow M/JM+I^{s+1}M \rightarrow 0.
\end{eqnarray*}
Since $M$ is Cohen-Macaulay, $\lambda(M/HM)=\lambda(M/JM)=e(I, M)$. By Proposition \ref{generic1}, $\lambda(M/HM+I^{s+1}M)\geq\lambda(M/JM+I^{s+1}M)$. Hence $\lambda(HM+I^{s+1}M/HM)\leq\lambda(JM+I^{s+1}M/JM).$ Hence, if $I^{s+1}M\subseteq JM$, then $\lambda(JM+I^{s+1}M/JM)=0$. Therefore, $\lambda(HM+I^{s+1}M/HM)=0$, which implies $I^{s+1}M\subseteq HM$.

To prove assertion (b), write $J=(x_1, \ldots, x_d)$ and $J_{d-1}=(x_1, \ldots, x_{d-1})M:_M I$, then
$$J_{d-1}\cap IM=(x_1, \ldots, x_{d-1})M$$
 and
 $$j(I, M)=\lambda(M/J_{d-1}+x_dM)=\lambda(M/J_{d-1}+IM)+\lambda(IM/JM).$$ Hence, $\lambda(IM/JM)$ does not depend on the choice of $J$. Since  $\lambda(IM/JM+I^{s+1}M)$ does not depend on $J$, one has $\lambda(JM+I^{s+1}M/JM)=\lambda(IM/JM)-\lambda(IM/JM+I^{s+1}M)$ does not depend on $J$. This implies $s_J(I, M)$ does not depend on  $J$.

Finally, observe that $H_{d-1}$ is a geometric $d-1$-residual intersection of $I$ on $M$ and
 $$H_{d-1}\cap I^2M=(y_1, \ldots, y_{d-1})IM.$$ Furthermore, $$j(I,M)=\lambda(\overline{M}/y_d\overline{M})=\lambda(I\overline{M}/y_dI\overline{M})=\lambda(IM/(y_1, \ldots, y_{d-1})M+I^2M)+\lambda(I^2M/HIM).
 $$
 Similarly, $j(I,M)=\lambda(IM/(x_1, \ldots, x_{d-1})M+I^2M)+\lambda(I^2M/JIM)$. By Proposition \ref{generic1}, one has $\lambda(IM/(x_1, \ldots, x_{d-1})M+I^2M)\leq \lambda(IM/(y_1, \ldots, y_{d-1})M+I^2M)$ and $\lambda(I^2M/JIM)\leq \lambda(I^2M/HIM)$. Hence $\lambda(I^2M/JIM)=\lambda(I^2M/HIM)$ does not depend on $J$. From the following exact sequences:
\begin{eqnarray*}
& 0\rightarrow HIM+I^{s+1}M/HIM \rightarrow I^2M/HIM \rightarrow I^2M/HIM+I^{s+1}M \rightarrow 0,\\
& 0\rightarrow JIM+I^{s+1}M/JIM \rightarrow I^2M/JIM \rightarrow I^2M/JIM+I^{s+1}M \rightarrow 0,
\end{eqnarray*}
one also has
$\lambda(HIM+I^{s+1}M/HIM)\leq\lambda(JIM+I^{s+1}M/JIM)$.
Let $s=s_J(I, M)$. If $s=0$, we are done. Otherwise, $I^{s+1}M\subseteq JM\cap I^2M=JIM$. Hence, $\lambda(JIM+I^{s+1}M/JIM)=0$. Therefore, $\lambda(HIM+I^{s+1}M/HIM)=0$, which in turn implies $I^{s+1}M\subseteq HIM\subseteq HM$.
\QED
\bigskip

The next several results generalize to $j$-stretched ideals of arbitrary height the corresponding results for stretched ideals proved by Rossi and Valla in \cite{RV3}. 
\begin{Lemma}\label{Properties}
Let $M$ be a Cohen-Macaulay module of dimension $d$ over a Noetherian local ring  with infinite residue field.
Let $I$ be an ideal that is $j$-stretched on $M$, but does not have minimal $j$-multiplicity on $M$. 
 Let $J=(x_1, \ldots, x_d)$ be a general minimal reduction of $I$ on $M$. Assume either $\lambda(M/IM)<\infty$ and $(x_1, \ldots, x_{d-1})M\cap I^2M\subseteq JIM$, or $\ell(I,M)=d$, $I$ satisfies $G_{d}$,  $AN^-_{d-2}$ on $M$, ${\rm depth}\,(M/IM)\geq 1$. Then
 \begin{itemize}
\item[(a)] $j(I,M)\geq \lambda(\overline{M}/I\overline{M})+h+1$, where $$h=h(I,\overline{M}),\, \overline{M}=M/J_{d-1},\, J_{d-1}=(x_1, \ldots, x_{d-1})M:_M I.$$
\item[(b)] For every $j\geq 1$, we have $I^{j+1}M=JI^jM+a^jbnR$, where $a,b\in I$ and $a,b$ are not in $J$ and $n\in M$.
\item[(c)] For every $j\geq 1$, we have $a^jbn\,\m \subseteq I^{j+2}M+JI^jM$.
\item[(d)] $IM=bnR+(JM:_M a)\cap IM$.
\end{itemize}
\end{Lemma}

\demo
(a)  By \cite[Lemma~2.1]{PX1},
$$
j(I,M)=e(I,\overline{M})= \lambda(I\overline{M}/I^2\overline{M})+\lambda(I^2\overline{M}/x_dI\overline{M}).
$$
By definition of $h$, this equals $\lambda(\overline{M}/I\overline{M})+h+\lambda(I^2\overline{M}/x_dI\overline{M})$.
Hence, to finish the proof of (a), we have to show that $\lambda(I^2\overline{M}/x_dI\overline{M})\geq 1$.
 However, this length is $0$ if and only if $I$ has minimal $j$-multiplicity on $M$. Since this option is ruled out, the assertion follows.

(b) We prove the assertion by induction on $j$. By  the proof of Theorem \ref{AN}, $1=\lambda(I^2M/JIM+I^3M)$. This implies that $I^2M=JIM+I^3M+abnR$ for some
$a,b\in I \setminus J$ and $n\in M$. By Nakayama's Lemma, $I^2M=JIM+abnR$, proving the statement in the case $j=1$. Now, for any $j\geq 1$,
assume by inductive hypothesis that $I^{j+1}M=JI^jM+a^jbnR$. We need to show that $I^{j+2}M=JI^{j+1}M+a^{j+1}bnR.$
This holds since $I^{j+2}M=I(JI^jM+a^jbnR)=JI^{j+1}M+a^j(bI)nR \subseteq JI^{j+1}M+a^j(JIM+abnR)=JI^{j+1}M+a^{j+1}bnR$, where this latter equality holds since $a^jJIM\subseteq JI^{j+1}M$.

The proofs of (c) and (d) are similar to the corresponding statements in Lemma \cite[2.4]{RV3}. We write them for sake of completeness.

(c) By induction on $j\geq 1$. The case $j=1$ follows from Theorem \ref{AN}.
Now, assume $j\geq 1$ and $(a^jbn)\m\subseteq I^{j+2}M+JI^jM$. Then,
$$(a^{j+1}bn)\m=a[(a^jbn)\m]\subseteq a[I^{j+2}M+JI^jM]\subseteq I^{j+3}M+JI^{j+1}M.$$

(d) Since $aIM\subseteq I^2M=JIM + abnR$, we have $IM\subseteq (JIM+abnR):_M a$. It is easy to check that $(JIM+abnR):_M a=(JIM:_M a) + bnR$. We then have $$IM\subseteq [(JIM:_M a) + bnR ]\cap IM=bnR+(JM:_M a)\cap IM\subseteq IM.$$
\QED
\bigskip
From the Specialization Lemma it can be seen that the numbers $\nu_j=\lambda(I^{j+1}M/JI^jM)$ are well-defined (that is, they do not depend on $J$(.)
\begin{Lemma}\label{non-inc}
Let $M$ be a Cohen-Macaulay module of dimension $d$ over a Noetherian local ring  with infinite residue field. Let $I$ be an ideal that is $j$-stretched on $M$.
Let $J=(x_1, \ldots, x_d)$ be a general minimal reduction of $I$ on $M$.
Assume either $\lambda(M/IM)<\infty$ and $(x_1, \ldots,x_{d-1})M\cap I^2M\subseteq JIM$, or $\ell(I,M)=d$, $I$ satisfies $G_{d}$,  $AN^-_{d-2}$ on $M$, ${\rm depth}\,(M/IM)\geq 1$.
Let $a,b,n$ be as in Lemma \ref{Properties}. Then
\begin{itemize}
\item[(a)] $\nu_j\leq \nu_{j-1}$ for every $j\geq 2;$
\item[(b)] set $\overline{\nu_j}=\lambda(I^{j+1}\overline{M}/JI^j\overline{M})$. Then $\overline{\nu_j}\leq \nu_j$ for every $j\geq 2$ and $\overline{\nu_1}=\nu_1$.
\end{itemize}
\end{Lemma}

\demo
(a) A proof similar to Lemma \ref{Properties} (b) shows that, for every $j\geq 2$, the following map, given by multiplication by $a$, is actually an epimorphism
$$I^{j}M/JI^{j-1}M \stackrel{\cdot a}{\longrightarrow} I^{j+1}M/JI^{j}M\longrightarrow 0.$$
Therefore, we obtain the desired inequality between the lengths of these two modules.

(b) For any $j$ we have the natural epimorphism
$$I^{j+1}M/JI^jM\longrightarrow I^{j+1}\overline{M}/JI^j\overline{M} \longrightarrow 0,$$
inducing the inequality $\overline{\nu_j}\leq \nu_j$.
The second part of the statement follows from the equalities
$$\begin{array}{lll}
\overline{\nu_1}&=\lambda(I^2\overline{M}/JI\overline{M}) &=\lambda(I^2M/JIM+ [((x_1,\dots,x_{d-1})M:_M I^{\infty})\cap I^2M])\\
                &=\lambda(I^2M/JIM)                       &=\nu_1
\end{array}$$
where  $[(x_1,\dots,x_{d-1})M:_M I^{\infty}]\cap I^2M\subseteq JIM$ (by assumption and \cite[Lemma~3.2.(f)]{PX1}).
\QED
\bigskip

Let $M$ be a finite module of dimension $d$ over a Noetherian local ring  with infinite residue field and  $I$  an ideal  with $\ell(I,M)=d$.
Let $J=(x_1, \ldots, x_d)$ be a general minimal reduction of $I$ on $M$. As before, 
 write $J_{d-1}=(x_1,\ldots,x_{d-1})M:_M I^{\infty}$ and $\overline{M}=M/J_{d-1}$. 
 By \cite{PX2}, the Hilbert function $H_{I,\,\overline{M}}(j)=\lambda_{R}(I^j\overline{M}/I^{j+1}\overline{M})$  
 does not depend on $J$. In particular, it is well-defined  the integer   
 $h(I,\,M):=\lambda(I\overline{M}/I^2\overline{M})-\lambda(\overline{M}/I\overline{M})$ which is dubbed the {\it embedding codimension} of $I$ on $M$. Moreover, it is easily checked that
$$j(I, M)=e(I, \overline{M})=\lambda(\overline{M}/I\overline{M})+h(I, M)+K-1, 
\mbox{ where } K-1=\lambda(I^2\overline{M}/x_dI\overline{M}).$$

The following corollary provides information on $K$.

\begin{Corollary}\label{K}
Let $M$ be a Cohen-Macaulay module of dimension $d$ over a Noetherian local ring  with infinite residue field. Let $I$ be an $R$-ideal that is $j$-stretched on $M$, but does not have minimal $j$-multiplicity. Let $J$ be a general minimal reduction of $I$ on $M$.
  Assume either $\lambda(M/IM)<\infty$ and $(x_1, \ldots, x_{d-1})M\cap I^2M\subseteq JIM$, or $\ell(I,M)=d$,
  $I$ satisfies $G_{d}$,  $AN^-_{d-2}$ on $M$, ${\rm depth}\,(M/IM)\geq 1$. Then
$$ K\geq 2, \quad\quad  \nu_1=K-1,\quad\quad I^{K}M\nsubseteq JM,\quad\quad I^{K+1}M\subseteq JM. 
$$
\end{Corollary}

\demo
By the proof of  Lemma \ref{non-inc} (b) we have $K-1=\lambda(I^2M/JIM)=\nu_1$.
By Lemma \ref{Properties} (a) we have $K-1\geq 1$, whence $K\geq 2$.
By Lemma \ref{Properties} (b) one obtains that
$$ P_{I,\,\overline{M}/x_d\overline{M}}=\lambda(\overline{M}/I\overline{M})+hz+z^2+
 \cdots+z^{j(I,M)-h+1-\lambda(\overline{M}/I\overline{M})}$$
Therefore, $K$ is the least positive integer with
$$I^{K+1}M\subseteq [((x_1, \ldots, x_{d-1})M:_M I)\cap I^{K+1}M]+I^{K+2}M+JM\subseteq I^{K+2}M+JM.$$ By Nakayama's Lemma, $K$ is the least positive integer with $I^{K+1}M\subseteq JM$.
\QED
\bigskip

The next result is the last ingredient that we need to characterize the Cohen-Macaulayness of ${\rm gr}_{I}(M)$ when $I$ is $j$-stretched on $M$. It shows that the inclusion  $I^{K+1}M\subseteq JI^jM$ is equivalent to some Valabrega-Valla equalities for small powers of $I$. More precisely,
\begin{Proposition}\label{VV}
Let $M$ be a Cohen-Macaulay module of dimension $d$ over a Noetherian local ring with infinite residue field. Let $I$ be an ideal that is $j$-stretched on $M$. Let $J=(x_1, \ldots,x_d)$ be a general minimal reduction of $I$ on $M$.  Assume either $\lambda(M/IM)<\infty$ and $(x_1,\ldots,x_{d-1})M\cap I^2M\subseteq JIM$ or $\ell(I,M)=d$, $I$ satisfies $G_{d}$,  $AN^-_{d-2}$ on $M$, and ${\rm depth}\,(M/IM)\geq 1$. Let $K$ be as above. Then for  any $0\leq j\leq K$, we have:
 \begin{itemize}
\item[(a)] $JM\cap I^{j+1}M=JI^jM+a^KbnR$, where $a,b,n$ are as in Lemma \ref{Properties}$;$
\item[(b)]  $I^{K+1}M\subseteq JI^jM$ if and only if $JM\cap I^{n+1}M=JI^nM$ for every $n\leq j$.
\end{itemize}
\end{Proposition}

\demo 
(a) By Lemma \ref{Properties} (b), $JM\cap I^{j+1}M=JM\cap(JI^{j}M+a^jbnR)=JI^jM + a^jbnR \cap JM$.
Hence, it is enough to prove that $a^jbnR \cap JM=a^KbnR$. Since by construction $a^KbnR \subseteq I^{K+1}M\subseteq JM$, we have
$a^KbnR \subseteq a^jbnR \cap JM$, proving one inclusion. To prove the other inclusion we use descending induction
on $j\leq K$. 
For $j=K$ the statement follows from Lemma \ref{Properties} (b). Now, assume $a^jbnR \cap JM\subseteq JI^jM+ a^KbnR$ for some $j\leq K$. We show that $$a^{j-1}bnR \cap JM\subseteq JI^{j-1}M+ a^KbnR.$$
From the discussion before Proposition \ref{VV} it follows that that $I^{j}M\not\subseteq JM$, hence
$$a^{j-1}bnR \cap JM\subseteq a^{j-1}bn\m \cap JM\subseteq(I^{j+1}M + JI^{j-1}M) \cap JM=(I^{j+1}M \cap JM ) + JI^{j-1}M.$$
By inductive hypothesis, this equals $JI^jM+a^KbnR + JI^{j-1}M=JI^{j-1}M+a^KbnR$.

The proof of assertion (b) is similar to the one of \cite[Lemma~2.5.(ii)]{RV3}.
\QED
\bigskip

\section{Cohen-Macaulayness and almost Cohen-Macaulayness of ${\rm gr}_I(R)$}

In this section we employ the results of the previous two sections to prove the two main results of this paper on the associated graded rings of $j$-stretched ideals. In Theorem \ref{CM} we characterize the Cohen-Macaulayness of the associated graded rings of $j$-stretched ideals in terms of the reduction number and index of nilpotency of $I$. In Theorem \ref{2} we prove a generalized version of Sally's conjecture for $j$-stretched ideals.
Our work is inspired by Rossi-Valla and Polini-Xie.
\medskip

The next theorem is the first main result of this section. It characterizes the $j$-stretched ideals $I$ such that ${\rm gr}_I(M)$ is Cohen-Macaulay. It generalizes widely the main results of \cite{S4}, \cite{RV3} and \cite{PX1}.

\begin{Theorem}\label{CM}
Let $M$ be a Cohen-Macaulay module of dimension $d$ over a Noetherian local ring with infinite residue field. Let $I$ be an ideal that is $j$-stretched on $M$. Let $J=(x_1, \ldots, x_d)$ be a general minimal reduction of $I$ on $M$. Assume either $\lambda(M/IM)<\infty$ and $(x_1,\ldots,x_{d-1})M\cap I^2M\subseteq JIM$ or $\ell(I,M)=d$, $I$ satisfies $G_{d}$,  $AN^-_{d-2}$ on $M$, and  ${\rm depth}\,(M/IM)\geq 1$.
 Let $K$ be as above. Then, the following statements are equivalent:
 \begin{itemize}
\item[(a)]   $T={\rm gr}_I(M)$ is Cohen-Macaulay$;$
\item[(b)]   $r_J(I,M)=s_J(I,M);$
\end{itemize}

Additionally, if $R$ is equicharacteristic, then $($a$)$ and $($b$)$ are also equivalent to
\begin{itemize}
\item[(c)]   $I^{K+1}M=HI^KM$ for some minimal reduction $H$ of $I$ on $M$.
\end{itemize}

\end{Theorem}

\demo
We first prove a couple of claims:
\begin{Claim}\label{C1}
 We have that $JM\cap I^{j+1}M=JI^jM$ for every $j\geq 0$ if and only if $I^{K+1}M=JI^KM$.
\end{Claim}
The forward direction is straightforward. Conversely, assume $I^{K+1}M=JI^KM$. By Proposition \ref{VV}, for every $j\leq K$, $JM\cap I^{j+1}M=JI^jM$. If $j\geq K$, then $I^{j+1}M=I^{j-K}I^{K+1}M=I^{j-K}JI^KM=JI^jM$ and we obtain $I^{j+1}M\cap JM=JI^jM$.

\begin{Claim}\label{C2}
Set $g= {\rm grade} (I, M)$. Let $x_1^{*}, \ldots, x_d^{*}$ be the initial forms of $x_1,\ldots,x_d$ in ${\rm gr}_I(R)$.
If $I^{K+1}M=JI^KM$, then $x_1^{*}, \ldots, x_g^{*}$ form a ${\rm gr}_I(M)$-regular sequence.
\end{Claim}

Since $x_1,\dots,x_g$ are general elements in $I$ and $g={\rm grade}(I,M)$, then $x_1,\dots,x_g$ form a regular sequence on $M$. By Valabrega-Valla criterion (see for instance \cite[Proposition~2.6]{VV} or \cite[Theorem~1.1]{RV}), Claim \ref{C2} is proved if we can show that $(x_1,\ldots,x_g)M\cap I^jM=(x_1,\ldots,x_g) I^{j-1}M$ for every $j\geq 1$.
The case $\lambda(M/IM)<\infty$ follows from \cite[2.6]{RV3}. So we may assume ${\rm dim}\,M/IM>0$.
We use induction on $j$ to prove $(x_1,\ldots,x_i)M\cap I^jM=(x_1,\ldots,x_i) I^{j-1}M$ for every $j\geq 1$ and $0\leq i\leq d$. This is clear if $j=1$. So, let us assume that $j\geq 2$ and the equality holds for $j-1$. We prove it by descending induction on $i\leq d$. Since $I$ is $j$-stretched on $M$ with $I^{K+1}M=JI^KM$, by Claim \ref{C1}, $JM\cap I^{j}M=JI^{j-1}M$, proving the case $i=d$. Now, assume by induction that $(x_1,\ldots,x_{i+1})M\cap I^{j}M=(x_1,\ldots,x_{i+1})I^{j-1}M$. Then,
\begin{eqnarray*}
&&(x_1, \ldots, x_{i})M\cap I^{j}M \\
&=&(x_1, \ldots, x_i)M\cap (x_1,\ldots,x_{i+1})I^{j-1}M \hspace{3.3cm}  
\mbox{ by induction on i}\\
&=&(x_1,\ldots,x_i)M\cap ((x_1,\ldots,x_{i})I^{j-1}M+x_{i+1}I^{j-1}M)\\
&=& (x_1,\ldots,x_{i})I^{j-1}M+(x_1,\ldots,x_i)M\cap x_{i+1}I^{j-1}M\\
&=& (x_1,\ldots,x_{i})I^{j-1}M+x_{i+1}[((x_1,\ldots,x_{i})M:_Mx_{i+1})\cap I^{j-1}M]\\
&=&(x_1,\ldots,x_{i})I^{j-1}M+x_{i+1}[(x_1,\ldots,x_{i})M\cap I^{j-1}M] \ \hspace{.02cm} 
 \mbox{ by Lemma \cite{PX1} (a)  and (e)}\\
&= &(x_1,\ldots,x_{i})I^{j-1}M+x_{i+1}(x_1,\ldots,x_{i}) I^{j-2}M \hspace{2.3cm}  \mbox{ by induction on j}\\
&\subseteq &(x_1,\ldots,x_i)I^{j-1}M.
\end{eqnarray*}
With the obvious inclusion $(x_1,\ldots,x_i)I^{j-1}M\subseteq (x_1, \ldots, x_{i})M\cap I^{j}M$, this finishes the proof.

\medskip

We are now ready to prove the theorem.

(a) $\Longleftrightarrow$ (b). The proof is similar to \cite[Theorem~3.8]{PX1}. Set  $\delta(I,M)=d-g$. We prove the equivalence of (a) and (b) by induction on $\delta$.  If  $\delta=0$, the assertion follows because we proved in Claim \ref{C2} that $x_1^{*}, \ldots, x_g^{*}$ form a ${\rm gr}_I(M)$-regular sequence.  Thus we may assume  that $\delta(I,M)\geq 1$ and the theorem holds for smaller values of $\delta(I,M)$. In particular, $d\geq g+1$. Again in both cases, since $x_1^{*}, \ldots, x_g^{*}$ form a ${\rm gr}_I(M)$-regular sequence, we may factor out $x_1, \ldots, x_g$ to assume $g=0$. Now $d=\delta(I,M) \geq 1$.
Set  $H_0=0:_M I$ and $\overline{M}=M/H_0$. By \cite[Lemma~3.2]{PX1} (b) and (c),  $\overline{M}$ is $d$-dimensional and Cohen-Macaulay. 
By \cite[Lemma~3.2]{PX1} (e), (a), (g) and (d), we have that
$IM\cap H_0=0,\, {\rm grade}\,(I, \overline{M})\geq 1$, $I$ still satisfies  $G_{d}$ and $AN^-_{d-2}$ on $\overline{M}$ and ${\rm depth}\,(\overline{M}/I\overline{M})\geq {\rm min}\{{\rm dim} \overline{M}/I\overline{M}, 1\}$.
Furthermore, ${\rm dim}\, M={\rm dim}\,\overline{M}=d$ and $IM\cap H_0=0$ imply that $\ell(I, \overline{M})=\ell(I, M)=d$. It is also easy to see that $I$ is stretched on $\overline{M}$ too.  Again by $IM\cap H_0=0$, there is a graded exact sequence
\begin{equation}\label{eq11}
0\rightarrow H_0\rightarrow {\rm gr}_I(M)\rightarrow {\rm gr}_{I}(\overline{M})\rightarrow 0.
\end{equation}

Since $\delta(I,\overline{M})=d-{\rm grade}\,(I, \overline{M})<d=\delta(I,M)$, by induction hypothesis  ${\rm depth}({\rm gr}_{I}(\overline{M}))\geq d$ if and only if ${I}^{K+1}\overline{M}=JI^K\overline{M}$.  Since $H_0\cap IM=0$, then also
$H_0\cap JI^KM=0$, and
\begin{equation}\label{powers} {I}^{K+1}\overline{M}/J{I}^K\overline{M}\cong I^{K+1}M/JI^KM.
\end{equation}
Finally, notice that ${\rm depth}(H_0)\geq d$ because $\overline{M}$ is $d$-dimensional Cohen-Macaulay. Therefore,
$$\begin{array}{lll}
{\rm gr}_I(M) \mbox{ is Cohen-Macaulay } & \Longleftrightarrow {\rm gr}_I(\overline{M}) \mbox{ is Cohen-Macaulay }& \mbox{(by (\ref{eq11}))}\vspace{0.1in} \\
& \Longleftrightarrow I^{K+1}\overline{M}=JI^K\overline{M} &\mbox{ by inductive hypothesis}\vspace{0.1in}\\
  & \Longleftrightarrow I^{K+1}{M}=JI^K{M} &\mbox{ by (\ref{powers}).}
\end{array}$$
By Corollary \ref{K}, we have $K=s_J(I,M)$, that is, $K$ is the index of nilpotency of $I$ on $M$ with respect to a general minimal reduction $J$.
Hence, $I^{K+1}M=JI^KM$ if and only if  $r_J(I,M)\leq s_J(I,M)$. This, in turn, is clearly equivalent to $r_J(I,M)=s_J(I,M)$ (one always have $r_J(I,M)\geq s_J(I,M)$). Therefore, we obtained
$${\rm gr}_I(M) \mbox{ is Cohen-Macaulay  if and only if } r_J(I,M)= s_J(I,M).$$

Now, assume $R$ is equicharacteristic. We prove that (b) $\Longleftrightarrow$ (c).  Clearly (b) implies (c).  To prove the converse notice that, for a general minimal reduction $J$, Lemma \ref{specializ} implies
$$
K=s_J(I, M)\leq r_J(I,M)\leq r_H(I,M).
$$
If (c) holds, then one has $r_H(I,M)=K$. In turn, this yields $K=s_J(I, M)= r_J(I,M)=K$.
\QED
\bigskip

Thanks to Theorem \ref{CM}, we recover one of the main results of a recent paper of Polini-Xie.
\begin{Corollary}$($\cite[Theorem~3.9]{PX1}$)$\label{PX1}
Let $M$ be a Cohen-Macaulay module of dimension $d$ over a Noetherian local ring $R$, let $I$ be an ideal with $\ell(I,M) = d$. Assume ${\rm depth}\, (M/IM) \geq {\rm min}\,\{{\rm dim}(M/IM),1\}$ and $I$ satisfies $G_d$ and $AN_{
d−2}^-$ on $M$. If $I$ has minimal $j$-multiplicity on $M$ then ${\rm gr}_I(M)$ is Cohen-Macaulay.
\end{Corollary}

\demo If $I$ has minimal $j$-multiplicity, it is easily seen that $r_J(I, M)=1$. Then, $K=1$ and a straightforward application of Theorem \ref{CM} concludes the proof.
\QED
\bigskip

\begin{Discussion}
Theorem \ref{CM} gives an effective condition to check the Cohen-Macaulayness of associated graded modules. One should remark that this is different from the criteria given by Johnson and Ulrich in \cite{JU}. In this latter paper they required  depth conditions on the powers $I^j$ for $1\leq j\leq r_J(I)$, whereas we do not have this assumption $($but 
we require the $j$-stretchedness of $I$$)$.
\end{Discussion}

Our second main result, Theorem \ref{2}, proves  conditions ensuring  the almost Cohen-Macaulayness  of ${\rm gr}_I(M)$. Since we plan to reduce to the two-dimensional case, we need the following theorem.

\begin{Theorem}\label{dim2}
Let $M$ be a $2$-dimensional Cohen-Macaulay module over a Noetherian local ring with infinite residue field. Let $I$ be an ideal that is $j$-stretched on $M$ and $J=(x_1, x_2)$  a general minimal reduction of $I$ on $M$.  Assume either $\lambda(M/IM)<\infty$ or $\ell(I,M)=2$, $I$ satisfies $G_{2}$,  $AN^-_{0}$ on $M$ and ${\rm depth}\,(M/IM)\geq 1$. Assume further there exists a positive integer $p$ such that
  \begin{itemize}
\item[(i)] $\lambda(JM\cap I^{j+1}M/JI^jM)=0$ for every $j\leq p-1;$
\item[(ii)]  $\lambda(I^{p+1}M/JI^pM)\leq 1$.
\end{itemize}
Then,
\begin{itemize}
\item[(a)] $x_1^*$ is regular on ${\rm gr}_I(M)^+;$
\item[(b)] ${\rm depth}\,( {\rm gr}_I(M) )\geq 1$ .
\end{itemize}
\end{Theorem}

\demo
 If ${\rm dim}\,M/IM=0$, then both claims follow from \cite[Theorem 4.4]{RV}. Thus, we may assume ${\rm depth}\,(M/IM)>0$.
Since $\lambda(I^{p+1}M/JI^pM)\leq 1$,
one has $I^{p+1}M=abnR+JI^pM$ for some $a\in I, b\in I^p, n\in M$ with $abn\not\in JI^pM$. For $j\geq p+1$, the multiplication by $a$ gives a surjective map from $I^{j+1}M/JI^{j}M$ to $I^{j+2}M/JI^{j+1}M$. Thus, $\lambda (I^{j+1}M/JI^{j}M)\leq 1$ for every $j\geq p$.

Notice that $x_1$ is regular on $IM$, since $(0:_M x_1)  \cap IM =0$ (Lemma \cite{PX1} (e)).  Thus, to prove that  $x_1^*$ is regular on ${\rm gr}_I(M)_+={\rm gr}_I(IM)$, we only need to show $x_1 IM\cap I^jIM=x_1 I^{j-1}IM$
for every $j\geq 1$ by \cite[Proposition~2.6]{VV} (see also \cite[Lemma~1.1]{RV}).
This is clear if $j=1$; hence we can  assume $j\geq 2$. Let $^{^{\tratto}}$ denote images
in $\overline{M}=M/x_1M$ and set  $s=r_{J}(I, I\overline{M})$. We claim that it is enough to show   $r_{J}(I, IM)=s$. Indeed, if $r_{J}(I, IM)=s$, then $x_1 IM\cap I^jIM=x_1IM\cap JI^{j-1}IM$ for any $j\geq 1$. This is clear  if $s\leq p-1$.
Assume $s>p-1\geq 1$. If $1\leq j\leq p-1$, then $x_1IM\cap I^{j}IM=x_1IM\cap JM\cap I^{j}IM=x_1IM\cap JI^{j-1}IM$.
If $p\leq j\leq s$, then $JI^{j-1}IM+x_1M\cap I^{j}IM=JI^{j-1}IM$. This follows from the following easy inequality of lengths
\vskip -.4cm
\begin{eqnarray*}
0&<&\lambda (I^j IM/ JI^{j-1}IM+(x_1M \cap I^{j} IM))\\
&=&\lambda(I^j IM/ JI^{j-1}IM) - \lambda(JI^{j-1}IM+(x_1M \cap I^{j}IM) / JI^{j-1}IM)\\
&=& 1 - \lambda(JI^{j-1}IM+(x_1M \cap I^{j}IM)/JI^{j-1}IM).
\end{eqnarray*}
On the other hand, if $j\geq s+1= r_{J}(I, IM)+1$, then  $I^jIM=JI^{j-1}IM$ and, therefore, $x_1 IM\cap I^jIM=x_1IM\cap JI^{j-1}IM$ for any $j\geq s+1$.

Now, applying  an argument similar to the one of  Theorem  \ref{CM}, we have  $x_1 IM\cap I^jIM=x_1 I^{j-1}IM$
for every $j\geq 1$.

To complete the proof of (a), we still need to to show that $r_{J}(I, IM)=s$. For this purpose we use the Ratliff-Rush filtration
$\widetilde{I^j} IM:=\cup_{t\geq 1}(I^{j+t}IM:_{IM} I^t)$ as it is done for ideals of definition (see \cite[Theorem~4.2]{RV}). As noticed earlier,  $x_1$ is regular on $IM$. Thus, for instance by  \cite[Lemma~3.1]{RV},
there exists an integer $n_0$ such that
$I^jIM=\widetilde{I^j} IM$ for $j\geq n_0$, and
\begin{equation}\label{RRF}
\widetilde{I^{j+1}} IM  :_{IM} x_1=\widetilde{I^{j}}IM \quad \mbox{for every } j\geq 0.
\end{equation}

As before, let $\overline{M}=M/x_1M$ and $^{^{\tratto}}$   denote images in $\overline{M}$.  There are two filtrations:
$$
\overline{\mathbb{M}}: I\overline{M}\supseteq  I^2\overline{M}\supseteq\ldots\supseteq I^{j} \overline{M}\supseteq \ldots
$$
and
$$
\overline{\mathbb{N}}:I\overline{M} \supseteq \overline{\widetilde{I} IM}:= \widetilde{I} \   \overline{IM}\supseteq\ldots\supseteq \widetilde{I^{j-1}} \overline{IM}\supseteq \ldots
$$
Notice that $\overline{\mathbb{M}}$ is an  $I$-adic filtration and $\overline{\mathbb{N}}$ is a good $I$-filtration on $I\overline{M}$ (see \cite[pag.~9]{RV} for the definition of  good filtration).
Furthermore, $I$ is an ideal of definition on $I\overline{M}$, i.e., $\lambda_{R}\,(I\overline{M}/  I^2\overline{M}) < \infty $. Indeed, $(x_1M :_M x_2) \cap IM=x_1M $  (see Lemma \cite{PX1} (e)), which in turn forces
$x_2 \in I$
to be regular on $I\overline{M}$. This  yields $\lambda_{R}(I\overline{M}/
I^2\overline{M})  \leq \lambda_{R}(I\overline{M}/ x_2 I\overline{M}) < \infty. $
Thus, we are  in the context of the filtrations as treated in \cite{RV}.
Since $ I^{j-1} \, I\overline{M}=\widetilde{I^{j-1}}I\overline{M}$ for $j\geq n_0$, the associated graded modules  ${\rm gr}_{\overline{\mathbb{N}}}(I\overline{M})$ and ${\rm gr}_{\overline{\mathbb{N}}}(I\overline{M})$
have the same Hilbert coefficients $e_0$ and $e_1$. Again, because there exists an element in  $I$
which is regular on $I\overline{M}$, by \cite[Lemmas~2.1 and 2.2]{RV}, we have
\begin{eqnarray*} \label{eq3}
&\sum_{j\geq 0}^{p-2} \lambda(I^{j+1}IM/JI^jIM)+s-(p-1)\\
&=\sum_{j\geq 0} \lambda(I^{j+1}I\overline{M}/ x_2 I^{j}I\overline{M}) =e_1(\overline{\mathbb {M}})=e_1(\overline{\mathbb {N}}) =
\sum_{j\geq 0} \lambda(\widetilde{I^{j+1}}I\overline{ M}/ x_2\widetilde{I^{j}} I\overline{ M}).
\end{eqnarray*}

 Observe that the first equality follows from the fact that, for $0\leq j \leq s$,
$$
\lambda(I^{j+1}M/JI^{j}M)=\lambda(I^{j+1}\overline{M}/JI^{j}\overline{M}).
$$
This holds because $\lambda(JM\cap I^{j+1}M/JI^{j}M)=0$ when  $0\leq j\leq p-1$. Hence,
\begin{eqnarray*}
&\lambda(I^{j+1}M/JI^{j}M)=\lambda(I^{j+1}M/JM\cap I^{j+1}M)=\lambda(I^{j+1}\overline{M}/J\overline{M}\cap I^{j+1}\overline{M})\\
&\leq \lambda(I^{j+1}\overline{M}/JI^j\overline{M})\leq \lambda(I^{j+1}M/JI^jM).
\end{eqnarray*}
On the other hand, if $p\leq j\leq s$, we have $0< \lambda(I^{j+1}\overline{M}/JI^j\overline{M}) \leq \lambda(I^{j+1}M/JI^jM)=1$. This proves that $\lambda(I^{j+1}\overline{M}/JI^j\overline{M}) = \lambda(I^{j+1}M/JI^jM)=1$ for $p\leq j\leq s=r_{J}(I, \overline{M})$
 and $\lambda(I^{j+1}\overline{M}/ JI^{j}\overline{M}) =0$ for $j \geq s+1$.

 \medskip

We now prove that $ \lambda(\widetilde{I^{j+1}}I\overline{ M}/ x_2\widetilde{I^{j}} I\overline{ M}) =  \lambda(\widetilde{I^{j+1}} IM/ J \widetilde{I^j} IM)$ for every $j\geq 0$.
Since
$$\widetilde{I^{j+1}}I\overline{ M}/ x_2\widetilde{I^{j}} I\overline{ M}\cong \widetilde{I^{j+1}} IM/(x_1M\cap \widetilde{I^{j+1}}IM+ x_2\widetilde{I^{j}} IM),
$$
we just need to show $x_1M\cap \widetilde{I^{j+1}}IM=x_1\widetilde{I^{j}} IM$.
We first prove that $x_1M\cap \widetilde{I}IM=x_1IM$.
Since $x_1M\cap \widetilde{I}IM\supseteq x_1IM$,
it suffices to show the equality locally at every associated prime ideal of $M/x_1IM$. By \cite[Lemma~3.2.(d)]{PX1}, every $\p\in {\rm Ass}(M/x_1IM)$ is not maximal. Hence, $IM_\p=\widetilde{I}M_\p=x_1 M_\p$, and, thus, $x_1 M_\p\cap \widetilde{I} IM_\p= \widetilde{I} IM_\p=x_1I M_\p$. Therefore, $x_1M\cap \widetilde{I}IM=x_1IM$. Now, for any $j\geq 1$, $x_1M\cap \widetilde{I^{j+1}}IM=x_1IM \cap \widetilde{I^{j+1}}IM=x_1(\widetilde{I^{j+1}}IM :_{IM} x_1)=x_1\widetilde{I^{j}} IM$.
Then,
\begin{equation}\label{eq4}
\sum_{j \ge 0}\lambda(\widetilde{I^{j+1}} IM/ J \widetilde{I^j} IM) =\sum_{j\geq 0}^{p-2} \lambda(I^{j+1}IM/JI^jIM)+s-(p-1).
\end{equation}

Let $W_J=\{n\geq 1\,\mid J\widetilde{I}^jIM\cap I^{j+1}IM=JI^{j}IM, \,0\leq j\leq n\}$. Then, $p-2\in W_J$. Hence, by \cite[Theorem~4.2]{RV},
$$
r(I, IM)\leq \sum_{j\geq 0}\lambda(\widetilde{I^{j+1}} IM/ J \widetilde{I^j} IM)+p-1-\sum_{j=0}^{p-2}\lambda (I^{j+1}IM/JI^{j}IM)=s.
$$
\smallskip

Finally, assertion (b) follows from (a). Indeed,  by assumption,  ${\rm depth}\,(M/IM)>0$ and from the exact sequence
$$
0\rightarrow M/IM \rightarrow {\rm gr}_I(M) \rightarrow {\rm gr}_I(M)^+ \rightarrow 0,
$$
we conclude that
$${\rm depth} ({\rm gr}_I(M))\geq {\rm min} \{{\rm depth} \, M/IM , {\rm depth} ({\rm gr}_I(M)^+) \} \geq 1.$$
\QED
\bigskip

The following theorem gives  a sufficient condition for the almost Cohen-Macaulayness of ${\rm gr}_I(M)$. It provides  a  generalized version of Sally's conjecture for arbitrary ideals and widely extends \cite[Theorem~4.4]{RV}.
\begin{Theorem}\label{2}
Let $M$ be a Cohen-Macaulay module of dimension $d$ over a Noetherian local ring with infinite residue field. Let $I$ be an ideal that is $j$-stretched on $M$ and  $J=(x_1, \ldots, x_d)$  a general minimal reduction of $I$ on $M$.  Assume either $\lambda(M/IM)<\infty$  or $\ell(I,M)=d$, $I$ satisfies $G_{d}$,  $AN^-_{d-2}$ on $M$ and ${\rm depth}\,(M/IM)\geq 1$. Assume there exists a positive integer $p$ such that
  \begin{itemize}
\item[(a)] $\lambda(JM\cap I^{j+1}M/JI^jM)=0$ for every $j\leq p-1;$
\item[(b)]  $\lambda(I^{p+1}M/JI^pM)\leq 1$.
\end{itemize}
Then ${\rm depth}\,( {\rm gr}_I(M) )\geq d-1$, i.e., ${\rm gr}_I(M)$ is almost Cohen-Macaulay.
\end{Theorem}

\demo    We prove the
theorem by induction on $d$. The case  $d=2$ being proven in Theorem
\ref{dim2}. Let $d\geq 3$ and assume the theorem holds for $d-1$. We
first reduce to the case  ${\rm grade}\, (I, M)  \geq1$. If ${\rm
grade}\,(I, M)=0$,  let $H_0=0:_M I$. As in the proof of Theorem \ref{CM} (see also
\cite{PX1}), all assumptions still hold for the module $M/H_0$.
Furthermore  $IM/H_0=IM/H_0 \cap IM=IM$, ${\rm grade} \, (I,
M/H_0) \geq 1$ and
${\rm depth}\,( {\rm gr}_I(M) )\geq {\rm depth}\,( {\rm
gr}_{I}(M/H_0) )$. So we are reduced to the case where the ideal $I$
has at least one  $M$-regular element.  Thus,  $x_1$ is regular on $M$.

If  ${\rm dim}\,M/IM=0$ then the assertion  follows from
\cite[Theorem~4.4]{RV}. Thus we may assume ${\rm dim}\,M/IM>0$. Let
$^{^{\tratto}}$  denote images in $\overline{M}=M/x_1M$.   Observe
that $\overline{M}$ is a Cohen-Macaulay module of dimension $d-1$
and $\ell(I, \overline{M})=d-1$. Also $I$ satisfies $G_{d-1}$ and
$AN^-_{d-3}$ on $\overline{M}$ by Lemma 3.2 in  \cite{PX1}.
Furthermore, observe  $\overline{M}/I\overline{M} \cong M/IM $  thus
${\rm depth}\,(\overline{M}/\overline{I})={\rm depth}\,(M/IM)\geq
{\rm min}\{{\rm dim} \, M/IM, 1\}= \{{\rm dim} \,
\overline{M}/I\overline{M}, 1\}$. Clearly,  $I$ is
$j$-stretched on $\overline{M}$ with respect to $\overline{J}$. By induction hypothesis,
$$
{\rm
depth}\,( {\rm gr}_{I}(\overline{M}) )\geq d-2.
$$

Next, we prove that $x_1^*$ is regular on ${\rm gr}_I(M)$. Since $x_1$ is regular on $M$, by \cite[Proposition~2.6]{VV} (see also \cite[Lemma~1.1]{RV}), the claim follows if the intersections $x_1M\cap I^{j}M=x_1I^{j-1}M$  hold for every $j\geq 1$. This is clear if $j=1$. If $j=2$, since $x_1IM\subseteq x_1M\cap I^2M$, it suffices to show the equality locally at every prime ideal $\p\in {\rm Ass}(M/x_1IM)$. By  Lemma \cite{PX1}, ${\rm depth}\,(M/x_1IM)\geq 1.$ Thus, for every prime ideal  $\p\in {\rm Ass}(M/x_1IM)$, $\p$ is not the maximal ideal of $R$ and, hence, either $IM_{\p}=M_{\p}$ or $IM_{\p}=(x_1, \ldots, x_{d-1})M_{\p}$. Therefore, $(x_1,\ldots,x_{d-1})M_{\p}\cap I^2M_{\p}=(x_1, \ldots, x_{d-1})IM_{\p}$. We use descending induction on $i$ to prove $(x_1,\ldots, x_{i})M_{\p}\cap I^2M_{\p}=(x_1, \ldots, x_{i})IM_{\p}$ for every $1\leq i\leq d-1$. Assume $(x_1,\ldots,x_{i+1})M_{\p}\cap I^2M_{\p}=(x_1,\ldots,x_{i+1})IM_{\p}$. Then
\begin{eqnarray*}
(x_1, \ldots, x_{i})M_{\p}\cap I^2M_{\p}&=&(x_1, \ldots, x_i)M_{\p}\cap (x_1,\ldots,x_{i+1})IM_{\p}\\
&=&(x_1,\ldots,x_i)M_{\p}\cap ((x_1,\ldots,x_{i})IM_{\p}+x_{i+1}IM_{\p})\\
&=& (x_1,\ldots,x_{i})IM_{\p}+(x_1,\ldots,x_i)M_{\p}\cap x_{i+1}IM_{\p}\\
&=& (x_1,\ldots,x_{i})IM_{\p}+x_{i+1}[((x_1,\ldots,x_{i})M_{\p}:_{M_{\p}}x_{i+1})\cap IM_{\p}]\\
&=& (x_1,\ldots,x_{i})IM_{\p}+x_{i+1}(x_1,\ldots,x_{i})M_{\p}\subseteq (x_1,\ldots,x_i)IM_{\p}.
\end{eqnarray*}
When $j\geq 3$, we have $x_1M\cap I^{j}M=x_1M\cap I^2M \cap I^{j}M=x_1IM\cap I^{j-1}IM=x_1I^{j-2}IM=x_1I^{j-1}M$, since $x_1^*$ is regular on ${\rm gr}_I(M)^+$ by Theorem \ref{dim2}.  Finally, since ${\rm depth}({\rm gr}_{I}(\overline{M}))\geq d-2$ and $x_1^*$ is regular on ${\rm gr}_I(M)$, we have ${\rm depth}({\rm gr}_{I}(M))\geq d-1$.
\QED
\bigskip

Notice that, thanks to the Specialization Lemma, in the assumption of Theorem \ref{2}.(b), one could replace a general minimal reduction of $I$ on $M$ by {\it any} minimal reduction of $I$ on $M$.

We can now state a concrete sufficient condition for the almost Cohen-Macaulayness of ${\rm gr}_I(M)$. 
\begin{Corollary}\label{3}
Let $M$ be a Cohen-Macaulay module over a Noetherian local ring with infinite residue field and ${\rm dim}\,M=d$. Let $I$ be an ideal $j$-stretched on $M$, $J=(x_1, \ldots, x_d)$ a general minimal reduction of $I$ on $M$, and $K$ be as above. Assume either $\lambda(M/IM)<\infty$ and $(x_1,\ldots, x_{d-1})M\cap I^2M\subseteq JIM$ or $\ell(I,M)=d$, $I$ satisfies $G_{d}$,  $AN^-_{d-2}$ on $M$ and ${\rm depth}\,(M/IM)\geq 1$, then
 \begin{itemize}
\item[(a)]   $I^{K+1}\subseteq JI^{K-1}$ if and only if $\lambda(I^{K}/JI^{K-1})=1$.
\item[(b)]  If $I^{K+1}\subseteq JI^{K-1}$,  then ${\rm depth}({\rm gr}_I(M))\geq d-1.$
\end{itemize}
\end{Corollary}

\demo
The proof is similar to \cite[Proposition~3.1]{RV3}. One needs to apply Theorem \ref{2}.
\QED
\bigskip

We now recover the second main result of \cite{PX1} as a special case of Theorem \ref{2}. 
\begin{Corollary}$($\cite[Theorem~4.8]{PX2}$)$\label{PX2}
Let $M$ be a Cohen-Macaulay module of dimension $d$ over a Noetherian local ring $R$.
Let $I$ be an ideal with $\ell(I,M) = d$. Assume ${\rm depth}\, (M/IM) \geq {\rm min}\{{\rm dim}\,(M/IM),1\}$ and $I$ satisfies $G_d$ and $AN^-_{d-2}$ on $M$. If $I$ has almost minimal $j$-multiplicity on $M$, then ${\rm depth}\,({\rm gr}_I(M)) \geq d-1.$
\end{Corollary}

\demo
If $I$ has almost minimal $j$-multiplicity on $M$, then $K=2$. Since $I^2\nsubseteq J$ and $\lambda(I^2/IJ)=1$, one has $I^2\cap J=IJ$. Therefore,  $I^3\subseteq JI$. Now, Corollary \ref{3} finishes the proof.
\QED
\bigskip

In \cite{R1} and \cite{RV3} it was introduced the concept of type of an ideal $I$ with respect to a given minimal reduction $J$ of $I$. This was defined as $\tau(I)=\lambda((J:I)\cap I/J)$. This definition depends heavily on the minimal reduction $J$.
Hence, we introduce a slight variation of this concept, that works also with respect to modules. For a general minimal reduction $J$ of $I$ on $M$, we set $$\tau(I,M)=\lambda((JM:_MI)\cap IM/JM),$$ and we call it the {\it general Cohen-Macaulay type} of $I$ on $M$.
Once again, thanks to the Specialization Lemma we are able to prove that when  the $G_d$ property is present this number is well-defined (in the sense that it is constant for $J$ general). This is achieved in the next Lemma.%

\begin{Lemma}\label{type}
Let $M$ be a $d$-dimensional Cohen-Macaulay module over a Noetherian local ring  with infinite residue field. 
Let $I$ be an ideal. Assume $\ell(I,M)=d$ and $I$ satisfies the $G_d$ property on $M$.

Then, the number $\tau(I,M)$
 is independent of the general minimal reduction $J$.
\end{Lemma}

\demo This can be proved  by a similar argument as in the proof of the Specialization Lemma. 
\QED
\bigskip

In the same spirit of the definitions given in \cite{PX1}, we say that an ideal $I$ has {\em almost almost minimal} $j$-multiplicity on $M$, if $\lambda(I^2\overline{M}/x_dI\overline{M})\leq 2$, or, equivalently, if $K\leq 3$.

Our next goal is to employ Theorem \ref{2} to prove that $j$-stretched ideals having almost almost mininimal $j$-multiplicity and small general Cohen-Macaulay type give rise to almost Cohen-Macaulay associated graded modules. This provides a general, higher dimensional version of results of \cite{RV3}. The first step in this direction consists in proving that $j$-stretched ideals of small general type satisfy the inclusion  $I^{K+1}M\subseteq JI^2M$.

\begin{Theorem}\label{SmallType}
Let $M$ be a Cohen-Macaulay module of dimension $d$ over a Noetherian local ring with infinite residue field. Let $I$ be an ideal that is $j$-stretched on $M$. Let $J=(x_1, \ldots, x_d)$ be a general minimal reduction of $I$ on $M$.  Assume either $\lambda(M/IM)<\infty$ and $JM\cap I^2M=JIM$ or $\ell(I,M)=d$, $I$ satisfies $G_{d}$,  $AN^-_{d-2}$ on $M$ and ${\rm depth}\,(M/IM)\geq 1$. Let $K$ be as above. If $\tau(I,M)< h+1-\lambda(\overline{M}/I\overline{M})$, then
$$  \nu_2=K-2,\,\,\,  JM\cap I^3M=JI^2M. $$
In particular $I^{K+1}M\subseteq JI^2M$.
\end{Theorem}

\demo
Similar to the proof of \cite[Theorem~2.7]{RV3}.
\QED
\bigskip

We now show that the associated graded modules of $j$-stretched ideals having almost almost mininimal $j$-multiplicity and small general Cohen-Macaulay type are almost Cohen-Macaulay. Since the cases $K=1,2$ have been proved in \cite{PX1}, we only need to prove the case $K=3$.
This result generalizes several classical results, see for instance \cite{S3}, \cite{RV2}, \cite{RV3} and \cite{RV}.

\begin{Corollary}\label{almalm}
Let $M$ be a Cohen-Macaulay module of dimension $d$ over a Noetherian local ring with infinite residue field. Let $I$ be an ideal that is $j$-stretched on $M$. Let $J$ be a general minimal reduction of $I$ on $M$.  Assume either $\lambda(M/IM)<\infty$ and $JM\cap I^2M=JIM$ or $\ell(I,M)=d$, $I$ satisfies $G_{d}$,  $AN^-_{d-2}$ on $M$, and ${\rm depth}\,(M/IM)\geq 1$.

If $K=3$ and $\tau(I,M)< h+1-\lambda(\overline{M}/I\overline{M})$, then $${\rm depth}({\rm gr}_I(M))\geq d-1.$$
\end{Corollary}

\demo By Theorem \ref{SmallType}, we have that $I^4M\subseteq JI^2M$. Now, we invoke Theorem \ref{2} to conclude.
\QED


\bigskip

\begin{thebibliography}{99}


\bibitem{AM}{R. Achilles and M. Manaresi, Multiplicity for ideals of maximal analytic spread and intersection theory,
J. Math. Kyoto Univ. {\bf 33-4} (1993), 1029--1046. }






\bibitem{CPV}{A. Corso, C. Polini and M. Vaz Pinto, Sally modules
and associated graded rings, Comm. in Algebra {\bf 26}
(1998), 2689--2708.}


\bibitem{E}{J. Elias, On the depth of the tangent cone and the
growth of the Hilbert function, Trans. Amer. Math. Soc. {\bf 351}
(1999), 4027--4042.}




\bibitem{F}{L. Fouli, A study on the core of ideals,  Ph.D. thesis, Purdue University, 2006.}






\bibitem{HKU}{W. Heinzer, M-K. Kim and B. Ulrich, The Gorenstein and complete intersection properties of associated graded rings, J. Pure Appl. Algebra {\bf 201} (2005), 264--283.}

\bibitem{H}{S. Huckaba, On the associated graded rings having almost
maximal depth, Comm. Algebra {\bf 26} (1998), 967--976.}

\bibitem{Hu}{ C. Huneke, Hilbert functions and symbolic powers, Michigan Math. J. {\bf 34} (1987), 293--318.}








\bibitem{JU}{M. Johnson and B. Ulrich, Artin--Nagata properties and
Cohen--Macaulay associated graded rings, Compositio Math. {\bf
103} (1996), 7--29. }


\bibitem{NT}{D. V.  Nhi and N. V. Trung, Specialization of modules, Comm. Algebra  {\bf 27} (1999), 2959--2978.}

\bibitem{NU}{K. Nishida and B. Ulrich, Computing $j$-multiplicities, to appear in J. Pure Appl. Algebra.}



\bibitem{P1} {T. Puthenpurakal,  Ratliff-Rush filtration, regularity and depth of
higher associated graded modules, Part I, J. Pure Appl. Algebra
{\bf 208} (2007), 159--176.}


\bibitem{PX1}{C. Polini and Y. Xie, $j$-multiplicity and depth of associated graded modules,  submitted.}

\bibitem{PX2}{C. Polini and Y. Xie, Generalized Hilbert functions, preprint.}






\bibitem{R1}{ M. E. Rossi, Primary ideals with good associated graded
rings, J. Pure Appl. Algebra {\bf 145} (2000), 75--90.}


\bibitem{RV1}{ M. E. Rossi and G. Valla, A conjecture of J. Sally, Comm. in Algebra
{\bf 24 (13)} (1996),  4249--4261.}

\bibitem{RV2}{ M. E. Rossi and G. Valla, Cohen-Macaulay local rings of embedding dimension $e+d-3$, J. London Math. Soc. {\bf 80} (2000), 107--126.}



\bibitem{RV3}{ M. E. Rossi and G. Valla, Stretched $\m$-primary ideals,  Beitr�ge Algebra Geom. {\bf 42 } (2001), 103--122.}


\bibitem{RV}{M. E. Rossi and G. Valla, Hilbert functions of filtered modules, Lecture Notes of the Unione Matematica Italiana, {\bf 9.} Springer-Verlag, Berlin; UMI, Bologna, 2010.}



\bibitem{S1}{ J. Sally, On the associated graded ring of a local Cohen-Macaulay ring,
J. Math. Kyoto Univ. {\bf 17} (1977),   19--21.}

\bibitem{S3}{J. Sally, Stretched Gorenstein rings, J. London Math. Soc. {\bf 20} (1979), 19--26.}

\bibitem{S2}{ J. Sally, Tangent cones at Gorenstein singularities,
Compositio Math. {\bf 40} (1980), 167--175.}


\bibitem{S4}{J. Sally, Cohen-Macaulay local rings of embedding dimension $e+d-2$,
J. Algebra {\bf 83} (1983), 393--408.}


\bibitem{T2}{N. V. Trung, Constructive characterization of the reduction numbers,
 Compositio Math.  {\bf 137}  (2003),  99--113.}

\bibitem{U}{B. Ulrich, Artin-Nagata properties and reductions of ideals,
Contemp. Math. {\bf 159} (1994), 373--400. }


\bibitem{VV}{P. Valabrega and G. Valla, Form rings and regular sequences,
Nagoya Math. J. {\bf 72} (1978), 91--101.}

\bibitem{V}{W. V. Vasconcelos, Hilbert functions, analytic spread and Koszul homology, 
Contemp. Math. {\bf 159}, Amer. Math. Soc., Providence, 1994, 401--422.}

\bibitem{Vaz}{ M. T. R. Vaz Pinto, Hilbert functions and Sally modules, J. Algebra {\bf 192} (1997),
504–-523.}

\bibitem{W}{H. J. Wang, On Cohen-Macaulay local rings with embedding dimension
e+d-2, J. Algebra {\bf 190} (1997),  226--240.}


\bibitem{X}{ Y. Xie, Formulas for the multiplicity of graded algebras, to appear in Trans.  Amer. Math. Soc. }

\end{thebibliography}
\end{document}